\RequirePackage{fix-cm}
\documentclass[smallcondensed]{svjour3}     
\smartqed  
\usepackage{graphicx}
\usepackage[scale=0.73]{geometry}
\usepackage{amsmath}
\usepackage{amssymb}
\usepackage{epic}
\usepackage{graphicx}
\usepackage{subfigure}
\usepackage{multirow}
\usepackage{float}
\usepackage{color}
\usepackage{graphicx}
\usepackage{listings}

\usepackage{enumerate}
\usepackage{hyperref}

\usepackage{verbatim}

 \journalname{SCIENCE CHINA Mathematics}

\begin{document}

\title{Asymptotic Coefficients and Errors for Chebyshev Polynomial Approximations with Weak Endpoint Singularities: Effects of Different Bases
\footnote{Accepted by Science China Mathematics on May 23th, 2022.}
}

\titlerunning{Asymptotic Coefficients and Errors for Chebyshev Approximations}        

\author{Xiaolong Zhang        \and
        John P. Boyd 
}


\institute{Xiaolong Zhang \at
              MOE-LCSM, School of Mathematics and Statistics, \\
            Hunan Normal University, Changsha, 410081, China \\
              \email{xlzhang@hunnu.edu.cn}           
           \and
           John P. Boyd \at
             Department of Climate \& Space
             Sciences and Engineering, University of Michigan,  \\ 2455 Hayward Avenue, Ann Arbor MI 48109 \\
             \email{jpboyd@umich.edu}
}

\date{Received: date / Accepted: date}

\maketitle

\begin{abstract}
	
When solving differential equations by a spectral method, it is often convenient to shift from 
Chebyshev polynomials $T_{n}(x)$ with coefficients $a_{n}$ to  modified basis functions that incorporate the boundary conditions. For homogeneous Dirichlet boundary conditions, $u(\pm 1)=0$, popular choices include the ``Chebyshev difference basis", $\varsigma_{n}(x) \equiv T_{n+2}(x) - T_{n}(x)$ with coefficients here denoted $b_{n}$ and the
``quadratic-factor basis functions"  $\varrho_{n}(x) \equiv (1-x^{2}) T_{n}(x)$ with coefficients $c_{n}$. If $u(x)$ is weakly singular at the boundaries, then $a_{n}$ will decrease proportionally to $\mathcal{O}(A(n)/n^{\kappa})$ for some positive constant $\kappa$, where the $A(n)$ is a logarithm or a constant. We prove that the Chebyshev difference coefficients $b_{n}$ decrease more slowly by a factor of $1/n$ while the quadratic-factor coefficients $c_{n}$ decrease more slowly still as $\mathcal{O}(A(n)/n^{\kappa-2})$.

The error for the unconstrained Chebyshev series, truncated at degree $n=N$,  is $\mathcal{O}(|A(N)|/N^{\kappa})$ in the interior, but is worse by one power of $N$ in narrow boundary layers near each of the endpoints. Despite having nearly identical error \emph{norms}, the error in the Chebyshev basis is concentrated in boundary layers near both endpoints, whereas the error in the quadratic-factor and difference basis sets is nearly uniform oscillations over the entire interval in $x$.

Meanwhile, for Chebyshev polynomials and the quadratic-factor basis, the value of the derivatives at the endpoints is $\mathcal{O}(N^{2})$, but only $\mathcal{O}(N)$ for the difference basis.

Furthermore, we have given the asymptotic coefficients and rigorous error estimates of the approximations in these three bases, solved by the least squares methods. In this paper, we also find an interesting fact that: on the face of it, aliasing error is regarded as a bad thing, actually, the error norm associated with the downward curving spectral coefficients decreases even faster than the error norm of infinite truncation.


\end{abstract}


\section {Introduction}

The success of Chebyshev polynomial spectral methods in solving differential and integral 
equations is comprehensively cataloged in a variety of standard texts such as \cite{FoxParker68,HesthavenGottliebGottlieb07,MasonHandscomb2002,Snyder66,Trefethen19}, and a cornucopia of others including two by the second author, \cite{Boyd99z} and \cite{BoydBook3}.  There are, however, some areas of spectral methods where open questions remain and consensus has not been achieved. One is the best way to impose boundary conditions. Even if we narrow the focus to the ``basis recombination", which is to use basis functions that are linear combinations of Chebyshev polynomials such that each basis function individually and exactly satisfies homogeneous linear boundary conditions, multiple options. Weak endpoint singularities --- `` weak" in the sense that the spectral series converges --- are still a topic of active exploration. In this article, we analyze both issues and show that they are closely 
interrelated.

The standard Chebyshev coefficients of a function $u(x)$ are the coefficients $a_{n}$ in the series
\begin{eqnarray}~\label{Eqa}
	u(x) = \sum_{n=0}^{\infty} a_{n}  T_{n}(x).
\end{eqnarray}
If $u(x)$ has weak endpoint singularities, then its Chebyshev coefficients $a_{n}$ will asymptotically (as $n \to \infty$) decrease proportional to $1/n^{\kappa}$ for some positive constant $\kappa$, which is the ``algebraic order of convergence", perhaps modulo some slower-than-power functions of $n$ such as $\ln^{\vartheta}(n), \vartheta \in \mathbb{N}_{+}$. Here, ``weak" [singularity] means that $u(x)$ is continuous everywhere on the interval $x \in [-1, 1]$, but its first derivative or higher derivatives are singular. 

When a problem satisfies homogeneous Dirichlet boundary conditions $u(\pm 1) = 0$, it is often {\it{desirable}} to choose basis functions that satisfy the boundary conditions. Two possibilities are
\begin{equation}\label{Eqb} 
	u^{diff}(x) = \sum_{n=0}^{\infty } b_{n} \,\varsigma_n(x), \quad 	\varsigma_n(x) \equiv T_{n+2}(x)-T_n(x) .\qquad
	\mbox{[Difference Basis]}
\end{equation}
or
\begin{equation}~\label{Eqc}
u^{quad}(x)  =  \sum_{n=0}^{\infty} c_{n} \, \varrho_{n}(x),\quad  	\varrho_{n}(x)  \equiv  (1-x^2) T_n(x).  \qquad \mbox{[Quadratic-factor Basis]}
\end{equation}

This was dubbed ``basis recombination" in the book of the second author, which discusses this 
strategy and its alternatives on pgs. 112 to 114 of \cite{Boyd99z}. The alternatives are ``boundary-bordering", which is to replace collocation or Galerkin projection conditions by rows of the discretization matrix that explicitly enforce the boundary conditions, and ``penalty methods" \cite{HesthavenGottliebGottlieb07}.
Karageorghis discusses the relationship between basis recombination and boundary bordering for multidimensional problems in single and multiple domains \cite{Karageorghis93b}.

Why ``desirable"? The second author gave an answer for eigenvalue problems in \cite{BoydOP4}. 
Boundary-bordering for an eigenvalue problem gives a discretization matrix in which the rows that impose the boundary conditions are independent of the eigenvalue. This was sufficient to wreck EISPACK, the premier eigensolver of its day. Forty years later, library matrix eigensolvers are made of sterner stuff, but the boundary imposing rows are still bad for the condition number.

Heinrichs \cite{Heinrichs91} pointed out that if one constructs the recombined basis functions to be, say,  for symmetric functions for Dirichlet boundary conditions $T_{n+2}(x) - T_{n}(x)$ instead of $T_{2 n+2}(x) - T_{0}(x)$, the oscillations of the two Chebyshev polynomials of similar degree partially cancel, reducing the condition number of the discretization matrix. The improvement for a $k$-th order differential equation with a basis truncated to $N$ Chebyshev polynomials is a factor of $N^{k}$ reduction in the condition number, which is particularly significant for high order differential equations.
\footnote{{\bf Parenthetically, note that basis recombination is also very convenient when $N$ is \emph{small} and the discretized problem is solved by a computer algebra system; reducing the number of basis coefficients from $N$ to $(N-2)$ greatly reduces the complexity of the explicit, analytic answer \cite{BoydOP55}.}}

The widespread use of basis recombination is attested by texts like \cite{HesthavenGottliebGottlieb07} as well as by other literature \cite{EllisonJulienVasil21}.

Convergence theory for Chebyshev polynomial series has co-evolved with Chebyshev 
algorithms and applications \cite{Boyd99z,HesthavenGottliebGottlieb07}. The second author's review summarizes convergence theory up to 2009 \cite{BoydOP173}. More recent contributions include \cite{Kzaz00,LiuWangLi19,Trefethen19,WangH16,WangHY2021,XiangLiu20}. 
There is also an active literature on closely related problems such as Gaussian quadrature and 
Clenshaw-Curtis for functions with various types of singularities, which was not included in \cite{BoydOP173} such as \cite{Riess1972,Trefethen08,WangHY18,XiangBornemann2012}. It is impossible to review this in detail, but the sheer mass of theory shows that this vein of mathematics is still being actively mined.

Gaps in the existing theory are: how do basis recombination and interpolation alter the convergence rate? In this article, we fill in these gaps.

One unnoticed but significant aspect of spectral methods for problems with weak 
endpoint singularities is that all three expansions have coefficients decreasing as inverse powers 
of $n$ (or inverse powers of $n$ multiplied by a factor of logarithm function), but the exponents are \emph{different} for each of the three as expressed by the first theorem below. Indeed, there are also other differences among these three of the basis sets when truncated.

Our comparisons employ three different ways to calculate the coefficients in these basis sets.
\begin{enumerate}[(1).]
	\item Chebyshev inner product projection,
	 \begin{eqnarray}\label{ChebCoeffsFormula}
	a_{0}= \frac{1}{\pi} \int_{-1}^{1} \,  \frac{T_{0}(x) }{\sqrt{1-x^{2}} } \, u(x) \, dx, 	\quad a_{n} = \frac{2}{\pi} \int_{-1}^{1} \,  \frac{ T_{n}(x) }{\sqrt{1-x^{2}} } \, u(x) \,  dx, \qquad n > 0 ,
	\end{eqnarray}
	$b_{n}$ and $c_{n}$ are whatever they need, as expressed  by the difference equations given 
	below, to be consistent to all degrees with the infinite Chebyshev series as defined precisely 
	in Theorem~\ref{Thabc}.
	
	\item The infinite sums are truncated and $(N+1)$-point interpolation is applied.
	
	\item Least squares minimization of constraints applied at $M$ points where $M>(N+1)$.
\end{enumerate}



\begin{figure}
	\centerline{\includegraphics[scale=0.75]{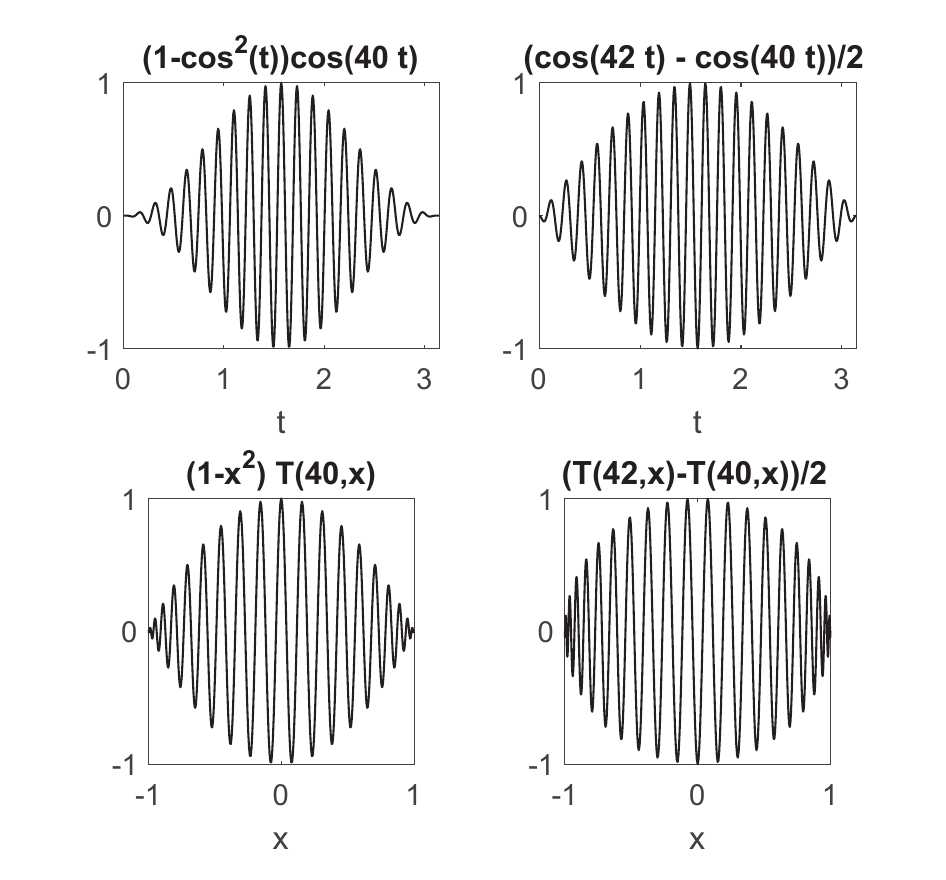}}
	\caption{The left two plots show a typical quadratic-factor basis function, plotted versus $t$ 
		at the top and $x=\cos(t)$ on the bottom. Right: Same but for a basis function which is 
		the difference of two Chebyshev polynomials, $\varsigma_{40}(x)= T_{42}(x) - T_{40}(x)$. The envelope of the bottom right curve is almost a circle of unit radius}
	\label{FigOP279_quadbasis_diffbasis_4plot}
\end{figure}

The least squares method yields a rectangular matrix problem. Interpolation is the limit
that the matrix is square, $M=N+1$, while the infinite series coefficients are the limit $M \rightarrow \infty$.

The effects on the \emph{errors} when each of the expansions is truncated 
after $N$ terms are subtle. These subtleties are explained in Sec. 3.

Fig.~\ref{FigOP279_quadbasis_diffbasis_4plot} compares basis functions. The qualitative 
resemblance is strong, which makes the behavioral differences all the more remarkable.

\section{Rates of decay of Chebyshev coefficients and  basis functions}~\label{SecRateCon}

In this section, we compare coefficients of the infinite series on each basis. Interpolation 
and least-squares with a finite number of quadrature points are reserved for later sections.

\begin{lemma}[Difference Equations for Infinite Series Coefficients]~\label{Thabc}
	Suppose a function $u(x)$ is zero at both endpoints but analytic everywhere on $ [-1, 1]$ except at the 
	endpoints where $u(x)$ is allowed to be weakly singular, ``weakly" in the sense that $v(x)=u(x)/(1-x^{2})$ is bounded at the endpoints. Let $u(x)$ have the three infinite series representations :  \eqref{Eqa}, \eqref{Eqb} and \eqref{Eqc}.
	Then  \begin{enumerate}[(i).]
		\item The $\varsigma_{n}(x)$ and $\varrho_{n}(x)$ are connected by the difference equation 
		and initial conditions
		\begin{equation}
		\begin{aligned}
			\varsigma_{0}(x) = T_{2}(x) - T_{0}(x) & =   - 2 \varrho_{0}(x) , \quad
			\varsigma_{1}(x) = T_{3}(x) - T_{1}(x)  =  -4 \varrho_{1}(x) , \quad   \\ 
			\varsigma_{n}(x) - \varsigma_{n-2}(x)  & =   - 4 \varrho_{n}(x), \quad n \ge 2.
	    \end{aligned}
		\end{equation}
		\item The condition $u^{diff}(x)=u(x)$ requires that the coefficients $b_{n}$ are connected to the
		$a_{n}$ by the difference equation:
		\begin{eqnarray}
			b_{0} & = & -a_{0}, \quad
			b_{1}  =  -a_{1},\quad
			b_{n-2} - b_{n}  =  a_{n}, \quad n\ge2, \label{Eq15}
		\end{eqnarray}
		which implies that
		\begin{eqnarray}
			b_{2n}  = - \sum_{j=0}^{n} a_{2j}, \quad
			b_{2n+1} = - \sum_{j=0}^{n} a_{2j+1}.
		\end{eqnarray}
		
		\item Similarly, $u^{quad}(x)=u(x)$ only if
		\begin{equation}\label{eqn:secondiff}
			\begin{aligned}
		&	\frac{1}{2} \, c_{0} \, - \, \frac{1}{4} c_{2} = a_{0}, \quad 
			\frac{1}{4}  c_{1} - \frac{1}{4} c_{3}  = a_{1},  \quad 
			-\frac{1}{2}  c_{0}-\frac{1}{4}  c_{4} + \frac{1}{2} c_{2} =  a_{2} ,\\
		&	-\frac{1}{4} \left( c_{n-2}+c_{n+2} \right)+ \frac{1}{2} c_{n}= a_{n},\quad (n  >  2).
			\end{aligned}
		\end{equation}

		\item The condition that $u^{quad}(x)=u^{diff}(x)$ demands that
		\begin{eqnarray}
			c_{0} - \frac{1}{2} c_{2}  =  - \, 2 \, b_{0} ,\quad 
			c_{n-2} - c_{n}  =   - 4 \, b_{n-2}, \qquad n \geq 3,
		\end{eqnarray}
		with the solution
		\begin{equation}~\label{Eqcnsol}
			c_{0}  =    -2 \sum_{j=0}^{\infty} b_{2j} , \quad 
			c_{1}  =   - 4 \sum_{j=0}^{\infty} b_{2j+1} , \quad 
			c_{2n}  =   - 4 \sum_{j=n}^{\infty} b_{2j}\,(n \ge 1) , \quad 
			c_{2n+1} = - 4 \sum_{j=n}^{\infty} b_{2j+1}. 
		\end{equation}
		Equivalently, using the infinite sums for $c_{0}$ and $c_{1}$, 
		the higher coefficients can be written as finite sums as
		\begin{eqnarray}
			c_{2n} & = &  2  c_{0} + 4 \sum_{j=0}^{n-1} b_{2j} \;( n \ge 1),\qquad
			c_{2n+1} =   c_{1} + 4 \sum_{j=0}^{n-1} b_{2j+1}. \qquad
		\end{eqnarray}
		
		\item Given $u(x) = u^{quad}(x)$, the coefficients $c_{n}$ can be defined without ambiguity as the Chebyshev coefficients 
		of an auxiliary function $v(x)$:
		\begin{eqnarray}\label{Eqn:vdef}
			v(x) \equiv \frac{u(x)}{1-x^{2}} = \sum_{n=0}^{\infty}\, c_{n} \, T_{n}(x),
		\end{eqnarray} 
		where  $c_{n}$ can be calculated by the formula \eqref{ChebCoeffsFormula}.
		
		\item If $u^{diff} (x) = u^{quad}(x)$, the relation of $b_{n}$ and $c_{n}$ is 
		\begin{eqnarray}
			b_{0} = \dfrac{c_{2} - 2 c_{0}}{4}  & = & \frac{1}{2 \pi} \int_{-1}^{1} \, \sqrt{ \frac{1}{1-x^{2}} }  \, v(x) \,
			\left( T_{2}(x) - T_{0}(x) \right) dx,\\
			b_{n} = \dfrac{c_{n+2} - c_{n}}{4}  & = & \frac{1}{2 \pi} \int_{-1}^{1} \, \sqrt{ \frac{1}{1-x^{2}} } \, v(x) \, \varsigma_{n}(x)
			\,   dx, \; n\in \mathbb{N}_{+}.
		\end{eqnarray}
	\end{enumerate}
\end{lemma}

\begin{proof}:
To show the first proposition, recall the Chebyshev identity \cite{Boyd99z,Snyder66}
\begin{equation*}
	T_{m}(x) T_{n}(x) = \frac{1}{2} \left( T_{m+n}(x) + T_{|m-n|}(x) \right).
\end{equation*}
It easily verify that the quadratic-factor basis can be written
\begin{eqnarray}~\label{Eqvarrho}
	\varrho_{0}(x) & = & \frac{1}{2} \left( T_{0}(x) - T_{2}(x) \right), \quad
	\varrho_{1}(x) =  \frac{1}{4} \left( T_{1}(x) - T_{3}(x) \right), \\
	\varrho_{n}(x) & = & (1-x^{2}) T_{n}(x) 
	=  - \frac{1}{4}  \left( T_{n+2}(x) + T_{n-2}(x)  \right) + \frac{1}{2}  \,  T_{n}(x) ,\quad n \geq 2. 
\end{eqnarray}
The difference of two difference basis functions is
\begin{eqnarray*}
	\varsigma_{n}(x) - \varsigma_{n-2}(x) & = & T_{n+2} (x)+ T_{n-2}(x) - 2\, T_{n}(x), \qquad n \ge 2,
\end{eqnarray*}
which is just $- 4 \varrho_{n}(x)$.

The second proposition follows from rewriting the series for $u^{diff}(x)$ as
\begin{eqnarray}
	u^{diff}(x) & = &  \sum_{n=0}^{\infty } b_{n} \left\{ T_{n+2}(x) \, - \, T_{n}(x) \right\} 
	 = - b_{0}  \, T_{0} (x) -  b_{1}  \,  T_{1} (x) +  \sum_{n=2}^{\infty } \, \left( b_{n-2} - b_{n} \right)  \, T_{n}(x) .
\end{eqnarray}
and comparing, term-by-term, with the standard Chebyshev series (\ref{Eqa}). The solution to the difference equation can be verified by direct substitution.

The reasoning to prove the third proposition, which is the second order difference equation for the $c_{n}$, is similar in that in the series for $u^{quad}(x)$, $\varrho_{n}(x)$ is replaced by its explicit expression in terms of Chebyshev polynomials, the sums are rearranged slightly as to extract the multiplier of $T_{n}(x)$, and this multiplier is equated with $a_{n}$:
\begin{eqnarray}
	u^{quad}(x)  =   \sum_{n=0}^{\infty} c_{n} \varrho_{n}(x)  
	&=&   \left(\,  \frac{1}{2} \, c_{0}  - \frac{1}{4} c_{2} \right) \, T_{0}(x) 
	+  \left( \, \frac{1}{4}  c_{1} - \frac{1}{4} \, c_{3} \right) T_{1}(x)
	+  \left(  \frac{1}{2} \, c_{2} - \,\frac{1}{2}  c_{0}\,  - \, \frac{1}{4} c_{4}   \right) \, T_{2}(x) 
	\nonumber \\   &&
	+ \sum_{n=3}^{\infty}\left( 
	-  \frac{1}{4}  c_{n-2} - \frac{1}{4} c_{n+2}   + \frac{1}{2} \, c_{n} \right)\,  T_{n}(x). \nonumber
\end{eqnarray}
Comparing this term-by-term with the Chebyshev series yields the difference equation. The
solution to the difference equation can again be verified by direct substitution.

The fourth proposition is demonstrated  by similarly rewriting the series for $u^{diff}(x)$ and $u^{quad}(x)$, substituting the expression for $\varrho_{n}(x)$ in terms of differences of $\varsigma_{n}(x)$, and then comparing the two series \cite{MasonHandscomb2002}.

Solving the recurrence is complicated because the lowest degree involves two `$c_{n}$'s. 
If we assume symbolic values for $c_{0}$ and $c_{1}$, we obtain the formal solution
\begin{eqnarray}
	c_{2n} & = &  2  c_{0} + 4 \sum_{j=0}^{n-1} b_{2j}\,(n \ge 1), \quad
	c_{2n+1} =   c_{1} + 4 \sum_{j=0}^{n-1} b_{2j+1}, \qquad
\end{eqnarray}
but this is not explicit without numerical values for $c_{0}$ and $c_{1}$.

On the other hand, if we truncate the infinite series so that $c_{N+1}=c_{N+2}=0$, then
\begin{eqnarray}
	c_{N}=-4 b_{N}, \qquad c_{N-1}=-4 b_{N-1}.
\end{eqnarray}
The recurrence can now be solved backwards to yield 
\begin{eqnarray}
	c_{0} & = &  - 2 \sum_{j=0}^{N_{e}} b_{2j}, \quad
	c_{2n}  =   - 4 \sum_{j=n}^{N_{e}} b_{2j}\,(n \ge 1) ,\quad
	c_{2n+1} = - 4 \sum_{j=n}^{N_{o}} n_{2j+1}, 
\end{eqnarray}
where $N_{e}=N_{o}=(N-1)/2$ if $N$ is odd and $N_{e}=N/2$ and $N_{o}=N/2-1$ if $N$ is even. The limit $N \to \infty$ yields the solution \eqref{Eqcnsol}.

The fifth proposition follows by dividing the series $u^{quad}(x)$ by $(1-x^{2})$ and then 
applying the usual integrals for Chebyshev coefficients.

Proposition Six follows from combining the difference relations connecting $b_{n}$ and $c_{n}$ (Proposition Four of this theorem) with the integrals for the $c_{n}$ proved as Proposition Five. $\blacksquare$
\end{proof}

Before analyzing the asymptotic decay rate of the Chebyshev coefficients of the infinite series for the functions with endpoint singularities. We shall give the exact representation of the Chebyshev coefficients for the function with an algebraic singularity. 

\begin{lemma}\label{lem:CoeffsAlg}\cite[(4.12)]{TuanElliott72}
	For the function $u(x)=(x+1)^{\varphi}$ with $\varphi > -\frac{1}{2} $ and $\varphi  \notin \mathbb{N}$ , the Chebyshev expansion coefficients are
	\begin{equation}
		a_{n} = \frac{(-1)^{n+1}}{\pi} \frac{\sin(\varphi\pi)}{2^{\varphi - 1}} \, \mathrm{B}(2\varphi+1, n-\varphi),\quad n\ge \varphi +1,
	\end{equation}
where $\mathrm{B}(x,y)$ denotes the Beta function.
\end{lemma}

The authors of \cite{LiuWangLi19} provided a detailed proof in the frame of fractional Sobolev-type spaces based on the generalized Gegenbauer functions of fractional degree (GGF-Fs). There is also other literature that discusses the consequences for orthogonal polynomial series to the function $u(x)$ with an algebraic singularity  \cite{BoydOP208,TuanElliott72,WangH16,WangHY2021,Xiang21,XiangLiu20}.

\begin{theorem}[Orders of Convergence for Coefficients of Infinite Series]~\label{Thm1}
	Suppose that $u(x)$ owns weak singularities at the endpoints as
	\begin{eqnarray}  \label{funsingular}
		u(x; \varphi, \vartheta) = g(x) (1-x^{2})^{\varphi}  \, \ln^{\vartheta}( 1 - x^{2}),\quad x\in[-1,1],
	\end{eqnarray} 
where $\varphi > \frac{1}{2}$, $\vartheta \in \mathbb{N}_{+}$,  $u(\pm1;\varphi,\vartheta)=\lim\limits_{x\to \pm 1}u(x;\varphi,\vartheta)=0$ and the function $g(x)$ is analytic everywhere on $x \in [-1, 1]$. Then the coefficients of three expansions \eqref{Eqa}, \eqref{Eqb} and \eqref{Eqc} respectively satisfy
	\begin{eqnarray}
		a_{n} & \sim & \frac{ A(n) } { n^{2 \varphi + 1} }, \\
		b_{n} & \sim & \frac{A(n)}{4 \, \varphi \, } \frac{1}{ n^{ 2 \varphi}},  \\
		c_{n} & \sim & - \frac{A(n) }{ (2 \, \varphi - 1) (2 \, \varphi)}   \frac{1}{ n^{ 2 \varphi - 1} }\label{eqn:quadraticCoeffs},
	\end{eqnarray}
	for $n\gg 1 $, where $A(n)$ varies more slowly than a power of $n$, such as a logarithm,
	or a constant. Specifically, $A(n)=\mathcal{O}\left(\ln^{\vartheta-1}(n)\right)$, when $\varphi \in \mathbb{N}$;  $A(n)=\mathcal{O}\left(\ln^{\vartheta}(n)\right)$ when $\varphi \notin\mathbb{N}$. Moreover, we have the expressions of $A(n)$:
	\begin{enumerate}[(i).]
		\item if $\vartheta=1$, one has
		\begin{eqnarray*}
			A(n)&=&-\frac{2\,\mathrm{\Gamma}(2\varphi+1) }{ \pi} \Big\{(-1)^{n}g(-1) + g(1)\Big\} \Big\{(\gamma_{1} - \ln 2 )\sin(\varphi\,\pi) + \pi \cos(\varphi\,\pi) \Big\} ,
		\end{eqnarray*}
	where $\gamma_{1}=2\,\psi_{0}(2\varphi +1) -\psi_{0}(n-\varphi)-\psi_{0}(n+\varphi+1)$ and  $\psi_{n}(x) (n\in \mathbb{N})$ is the polygamma function. 
	\item if $\vartheta=2$, one has
	\begin{equation*}
		\begin{aligned}
		A(n)=& -\frac{2\,\mathrm{\Gamma} (2\varphi +1 )}{\pi} \Big\{ (-1)^{n}g(-1) + g(1) \Big\} \left\{ \left(\gamma_{1}^{2} + \gamma_{2} -2\ln(2)\,\gamma_{1} + \ln^{2}(2) - \pi^2 \right) \sin(\varphi\,\pi) \right.\\
	  &	\left. + 2\pi \Big ( \gamma_{1} - \ln(2)\Big )\cos(\varphi\,\pi)\right \},
		\end{aligned}
	\end{equation*}
where $\gamma_{2}=4 \psi_{1} ( 2 \varphi + 1) + \psi_{1} ( n - \varphi ) - \psi_{1}( n + \varphi + 1 )$.

\item if $\vartheta\in \mathbb{N}_{+}$, one has the general formula 
\begin{equation*}
\begin{aligned}
\hspace{-0.2in}
A(n) =&	- \frac{2}{\pi} \,\gamma_{3}\, n^{2\varphi + 1} \sum_{k=0}^{\vartheta} \sum_{j=0}^{k} \binom{k}{j} \pi^{j} \sin\left( \frac{j}{2} \pi + \varphi\,\pi\right )\ln^{ k - j }  \left(\frac{1}{2} \right) \frac{d^{\vartheta -k }}{d\varphi^{\vartheta -k }} \mathrm{B}( 2 \varphi + 1, n - \varphi )\\
\sim  & -\frac{2\,\mathrm{\Gamma} (2\varphi +1 )}{\pi} \,\gamma_{3} \, \Big\{ \ln^{\vartheta} (n) \sin(\varphi \, \pi ) + \ln^{\vartheta - 1 } (n) \cos(\varphi \, \pi ) \Big \},
\end{aligned}
\end{equation*}
where
\begin{equation*}
	\gamma_{3} =  (-1)^{n}g(-1) + g(1).
\end{equation*}
\end{enumerate}

\end{theorem}

\begin{proof}
: The asymptotic behavior of the Chebyshev coefficients $a_{n}$ follows from a theorem of Elliott 
\cite{Elliott64}, but see also \cite{BoydOP50,ElliottTuan74,Kzaz00,LiuWangLi19,TuanElliott72,WangHY18,Xiang21,XiangLiu20}.  For simplicity, take $A(n)$ and $B(n)$ as constants below. Then for large $n$ and assuming power-law behavior for the $b_{n}$ with algebraic order of convergence $k$ with proportionality constant $B$, the difference equation (\ref{Eq15}) gives
\begin{eqnarray*}
	\frac{B}{n^{k}} \left( \frac{1}{ (1 - 2/n)^{k} } - 1 \right)  =  \frac{A}{n^{2 \varphi + 1}} .
\end{eqnarray*} 
For large $n$, $(1-2/n)^{-k} \approx 1 + 2 k / n + \mathcal{O}(k^{2}/n^{2})$ from
whence
\begin{eqnarray*}
	\frac{2 k B}{n^{k+1}} & = &  \frac{A}{n^{2 \varphi + 1}}, 
\end{eqnarray*}
from which it follows that $k = 2 \varphi$ as claimed and $B=A/(2k)= A/(4 \varphi)$.

%
%

To prove the third proposition, define $v(x)$ as before by $u(x) = (1 - x^{2}) \, v(x)$. 
The $c_{n}$ are the \emph{standard} Chebyshev polynomial coefficients of the 
\emph{modified}  function
\begin{eqnarray*}
	v = \sum_{n=0}^{\infty} c_{n}  T_{n}(x),
\end{eqnarray*}
for which the asymptotic behavior of the $c_{n}$ follows from Elliott's theorem \cite{Elliott64}.

An alternative proof that gives the relative proportionality constant is as follows. Earlier, we proved in Proposition 3 of Theorem~\ref{Thabc}  that
\begin{eqnarray}~\label{Eqzipbis}
	-\frac{1}{4} \left( c_{n-2}+c_{n+2} \right)+ \frac{1}{2} c_{n}=a_{n},\quad (n  \ge  2).
\end{eqnarray}
Assume that asymptotically for large degree $n$
\begin{equation*}
	c_{n} \sim \frac{C}{n^{k}},\quad n\gg 1.
\end{equation*}
Substituting this into the second order difference equation \eqref{eqn:secondiff} gives
\begin{equation*}
	\frac{C}{2} \frac{1}{n^{k}} \left(-\frac{1}{2}\left(\frac{1}{(1-2/n)^{k}}+\frac{1}{(1+2/n)^{k}}\right)+1\right)=\frac{A}{n^{2\varphi + 1}}.
\end{equation*}
Then, it is not difficult to see that
\begin{eqnarray*}
	\frac{Ck(k+1)}{n^{k+2}}& = & \frac{A}{n^{2 \varphi + 1}}.
\end{eqnarray*}
Thus it follows that $k = 2 \varphi-1$ and $C=- A/\{ k (k+1) \}= - A/(4 \varphi -2) \, \varphi)$.
The proof is not substantially changed if $A$ is allowed to vary slowly with the degree,  such as logarithmically. 

Recently, Liu et.al. give the optimal decay rate of the Chebyshev expansion coefficients for this function $u(x;\varphi, \vartheta) = (1+x)^{\varphi} \ln^{\vartheta}(1+x)$ when $\vartheta=1$ \cite{LiuWangLi19}.  By the idea of this paper, we will prove the optimal estimates of $A(n)$ given above, for the function \eqref{funsingular} when $\vartheta =1,2$.  By the Eq. \eqref{ChebCoeffsFormula}, for $n>0$, the Chebyshev expansion coefficients are
\begin{equation*}
	\begin{aligned}
a_{n} =& \frac{2} {\pi} \int_{-1}^{1}  g(x) (1-x^2)^{\varphi} \ln^{\vartheta}( 1 - x^2 ) \,\frac{T_{n}(x)}{\sqrt{ 1 - x^2 }}dx\\ 
 =& \frac{2} {\pi} \int_{-1}^{1}  g(x) (1 - x^{2})^{\varphi} \left\{ \ln^{\vartheta}( 1 + x )+\ln^{\vartheta}(1-x)+\sum_{i=1}^{\vartheta-1}\binom{\vartheta}{i}\ln^{i}(1+x)\ln^{\vartheta-i}(1-x)\right\} \frac{T_{n}(x)}{\sqrt{ 1 - x^2 }}dx.\\ 
    \end{aligned} 
\end{equation*}
As we known, the coefficients are dominated by terms who own the worst singularities, that is, the terms whose lowest order derivatives are unbounded and increase highest at the corresponding singularities. Thus, for $\vartheta >1$ 
\begin{equation}
	\begin{aligned}
	a_{n}& \simeq  \frac{2} {\pi} \int_{-1}^{1}  g(x) (1 - x^{2})^{\varphi} \left\{ \ln^{\vartheta}( 1 + x ) + \ln^{\vartheta}(1-x)  \right\} \,\frac{T_{n}(x)}{\sqrt{ 1 - x^2 }}dx\qquad \\
	&=\frac{2}{\pi} \int_{-1}^{1} g(x)(1-x^2)^{\varphi}\ln^{\vartheta}(1+x)\frac{T_{n}(x)}{\sqrt{ 1 - x^2 }}dx + \frac{2}{\pi} \int_{-1}^{1} g(x)  (1-x^{2})^{\varphi} \ln^{\vartheta}(1-x)\frac{T_{n}(x)}{\sqrt{ 1 - x^2 }}dx.	
	\end{aligned}
\end{equation}
For convenience, we set 
\begin{equation*}
	a_{n,1} := \frac{2}{\pi} \int_{-1}^{1} g(x)(1-x^2)^{\varphi}\ln^{\vartheta}(1+x)\frac{T_{n}(x)}{\sqrt{ 1 - x^2 }}dx,\quad 
	 a_{n,2}:=\frac{2}{\pi} \int_{-1}^{1} g(x)  (1-x^{2})^{\varphi} \ln^{\vartheta}(1-x)\frac{T_{n}(x)}{\sqrt{ 1 - x^2 }}dx.
\end{equation*}
To obtain the asymptotic behavior of $a_{n,1}$, the dominant term of the integrand needs to be considered. Due to the function $g(x)$ is analytic on the interval $[-1,1]$, it can be written as Taylor series at $x=-1$. Thus, the dominant contribution comes from the integral 
\begin{equation}
	a_{n,1} \simeq  \frac{2^{\varphi +1}}{\pi} g(-1)\int_{-1}^{1} \,(1+x)^{\varphi}\ln^{\vartheta}(1+x)\frac{T_{n}(x)}{\sqrt{ 1 - x^2 }}dx.
\end{equation}
By the Lemma \ref{lem:CoeffsAlg}, and using the L'Hospital rule, yields
\begin{equation*}
	\begin{aligned}
	a_{n,1} \simeq & \frac{2^{\varphi +1}}{\pi}  g(-1) \int_{-1}^{1} \,(1+x)^{\varphi}
	 \left\{  \lim_{\varepsilon \to 0}  \frac{  \left( (1+x)^{\varepsilon} -\sum_{j=0}^{\vartheta-1} \frac{ \ln^{j} (1+x) }{j!} \varepsilon^{j} \right) }{\frac{ \varepsilon^{\vartheta}} {\vartheta!}}\right\}
	\frac{T_{n}(x)}{\sqrt{ 1 - x^2 }}dx \\
	= &(-1)^{n+1} \frac{2}{\pi} g(-1) \sum_{k=0}^{\vartheta} \sum_{j=0}^{k} \binom{k}{j} \pi^{j} \sin\left( \frac{j}{2} \pi + \varphi\,\pi\right )\ln^{ k - j }  \left(\frac{1}{2} \right) \frac{d^{\vartheta -k }}{d\varphi^{\vartheta -k }} \mathrm{B}( 2 \varphi + 1, n - \varphi ).
	\end{aligned}
\end{equation*}
Similarly, the $a_{n,2}$ is just by a constant factor $\frac{(-1)^{n}g(1)}{g(-1)}$ of $a_{n,1}$.

 Recall that if $y$ is large and $x$ is fixed, then   
\begin{equation}\label{BetaAsym}
	\mathrm{B}(x,y)\sim \mathrm{\Gamma}(x) \,y^{-x}.
\end{equation}
By induction method, and using \eqref{BetaAsym}, it leads to 
\begin{equation*}
	A(n) \sim -\frac{2\,\mathrm{\Gamma} (2\varphi +1 )}{\pi} \Big\{ (-1)^{n}g(-1) + g(1) \Big\} \Big\{ \ln^{\vartheta} (n) \sin(\varphi \, \pi ) + \ln^{\vartheta - 1 } (n) \cos(\varphi \, \pi ) \Big \}. 
\end{equation*}
Note that, using the above method,  for $\vartheta=1$  the exact Chebyshev coefficients $a_{n}$ can be obtained; for $\vartheta=2$ we can obtain the optimal estimate of $A(n)$, that is, the dominated terms can be exactly achieved. For a general $\vartheta\in\mathbb{N}$, the rough estimates of the Chebyshev coefficients can be found in \cite{Xiang21,Zhang2021}.
$\blacksquare$
\end{proof}

Here, if $\vartheta=0$, then $u(x;\varphi,\vartheta)$ is singular only if $\varphi$ is not an integer, as is mentioned before this theorem.  If not otherwise specified, $A(n)$ in the remainder of this paper denotes the expression of $A(n)$ given in Theorem \ref{Thm1}.

\begin{remark}
Here, the parameter $\varphi>\frac{1}{2}$ is required to make sure $2\varphi -1>0$ in (\ref{eqn:quadraticCoeffs}). In fact, when using the Chebyshev basis to approximate the function (\ref{funsingular}), it is only required $\varphi> -\frac{1}{2}$.
\end{remark}

\begin{remark}
	Because the natural logarithm function $\ln(n)$ increases very slowly as $n$ increases, the differences in plots between $\ln(n)/n^{\kappa}$ and $1/n^{\kappa}$ are subtle in numerical experiments. It is very easy to believe that $A(n)$ is always a constant for all $\varphi>-\frac{1}{2}$. However, as demonstrated in this theorem, the $A(n)$ is a constant only when $\vartheta=1$ and $\varphi\in\mathbb{N}$.
\end{remark}
\begin{figure}
	\centerline{\includegraphics[width=7in,height=3.5in]{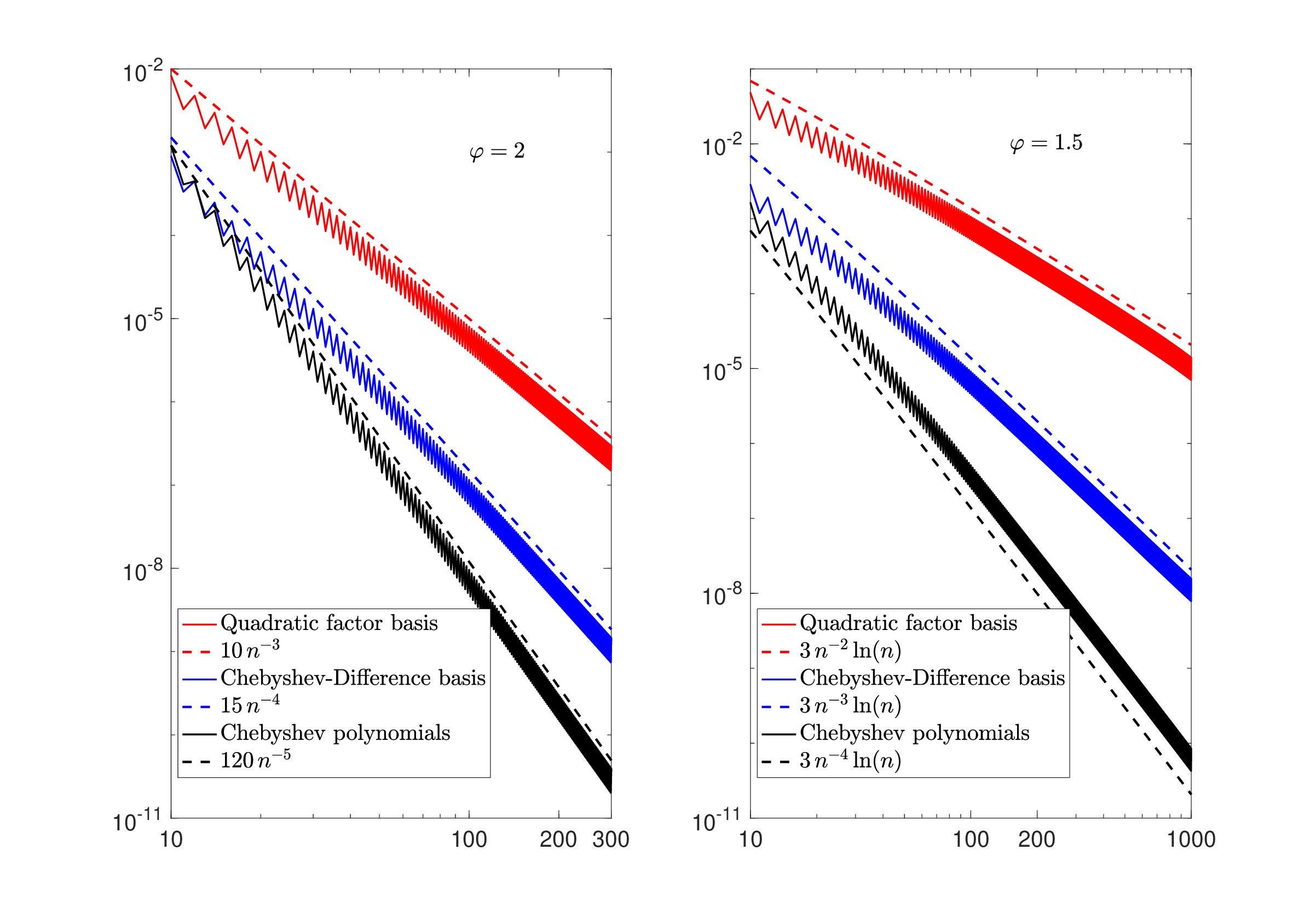}}
	\caption{Coefficients of the infinite series of the function $u(x) \equiv  (1+\frac{x}{2})(1-x^{2})^{\varphi} \, \ln(1 - x^{2})$ in three different bases with $\varphi =2$ (Left), $\varphi=1.5$ (Right).}
	\label{FigOP279_Coeffs_InfiniteSeries}
\end{figure}

Fig.~\ref{FigOP279_Coeffs_InfiniteSeries} confirms coefficients' law in Theorem \ref{Thm1}. Care must be exercised in interpreting this theorem. It applies when the $a_{n}$ obeys an inverse 
power law, as is true of the exact Chebyshev coefficients of the infinite series. We later compute a variety of finite approximations to $u(x)$ and these, when represented in the Chebyshev basis, 
do not automatically have the inverse power-law behavior of the coefficients $a_{n}$.

\section{Errors in Truncating Infinite Series}

Suppose we truncate each of the three series to a polynomial of degree $N$
\begin{eqnarray}
	u_{N}(x) & = & \sum_{n=0}^{N} a_{n}  T_{n}(x), ~\label{eqn:trunCheb}\\
	u_{N}^{diff}(x) & = &\sum_{n=0}^{N-2} b_{n} \left\{ T_{n+2}(x) \, - \, T_{n}(x) \right\}, \qquad
	\mbox{[Difference Basis]}\label{eqn:trunDiff}  \\
	u_{N}^{quad}(x)& = &\sum_{n=0}^{N-2} c_{n} (1-x^{2}) T_{n}(x) . \qquad \mbox{[Quadratic-factor Basis]} \label{eqn:trunQuad}
\end{eqnarray}
We have previously described the behavior of the \emph{coefficients} $a_{n}$, $b_{n}$ and $c_{n}$, but here the question is: what are the \emph{errors} in these truncations?

For the class of functions \eqref{funsingular},  Theorem~\ref{Thm1} shows that the Chebyshev coefficients fall as $\mathcal{O}(A(n)/n^{2 \varphi+1})$ 
while the quadratic-factor basis coefficients $c_{n}$ decrease as $\mathcal{O}(A(n)/n^{2 \varphi - 1})$. A well-known theorem asserts that truncation error in a Chebyshev series is bounded by the 
sum of the absolute values of all the neglected \emph{terms}; because $|T_{n}(x)| \leq 1$ on $x \in [-1, 1]$, the bound is also the sum of the absolute values of all the neglected \emph{coefficients}. One might suppose that the error in the $L_{\infty}$ norm when the series is truncated at $n=N$ is the magnitude of the largest omitted coefficient, but in fact, the series error is worse by $\mathcal{O}(1/N)$ than the rate of convergence of the Chebyshev coefficients. Near the endpoints, the terms are all of the same sign or asymptotically strictly alternating. The order of convergence of the error then comes from the asymptotic sum approximation (\ref{EqBasz}) below.

\begin{figure}[h]
	\centerline{\hspace{0.2in}\includegraphics[scale=0.6]{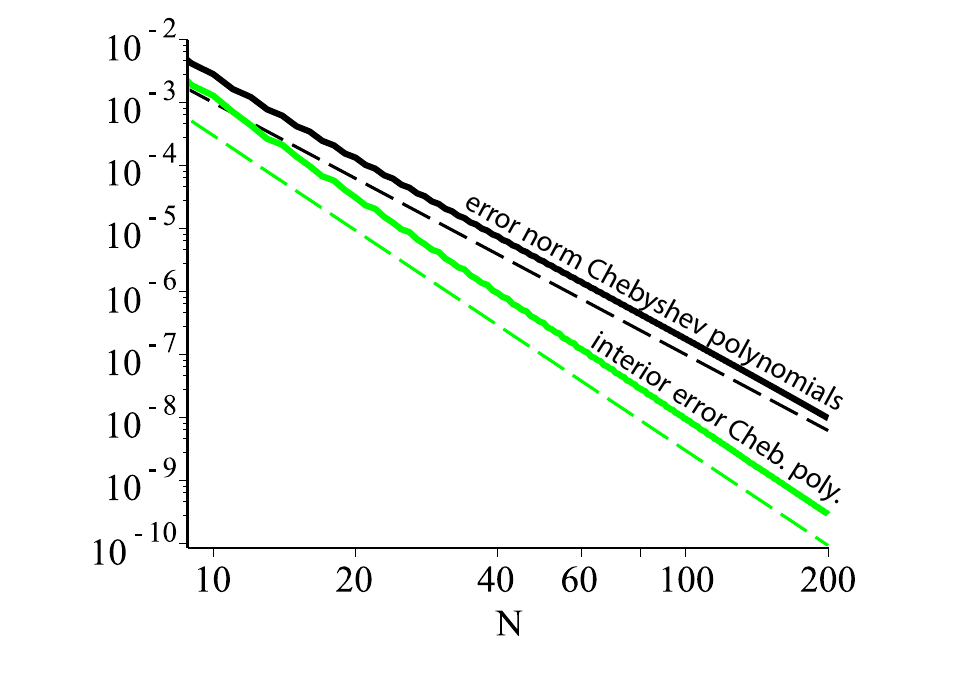},\hspace{-0.7in}\includegraphics[scale=0.5]{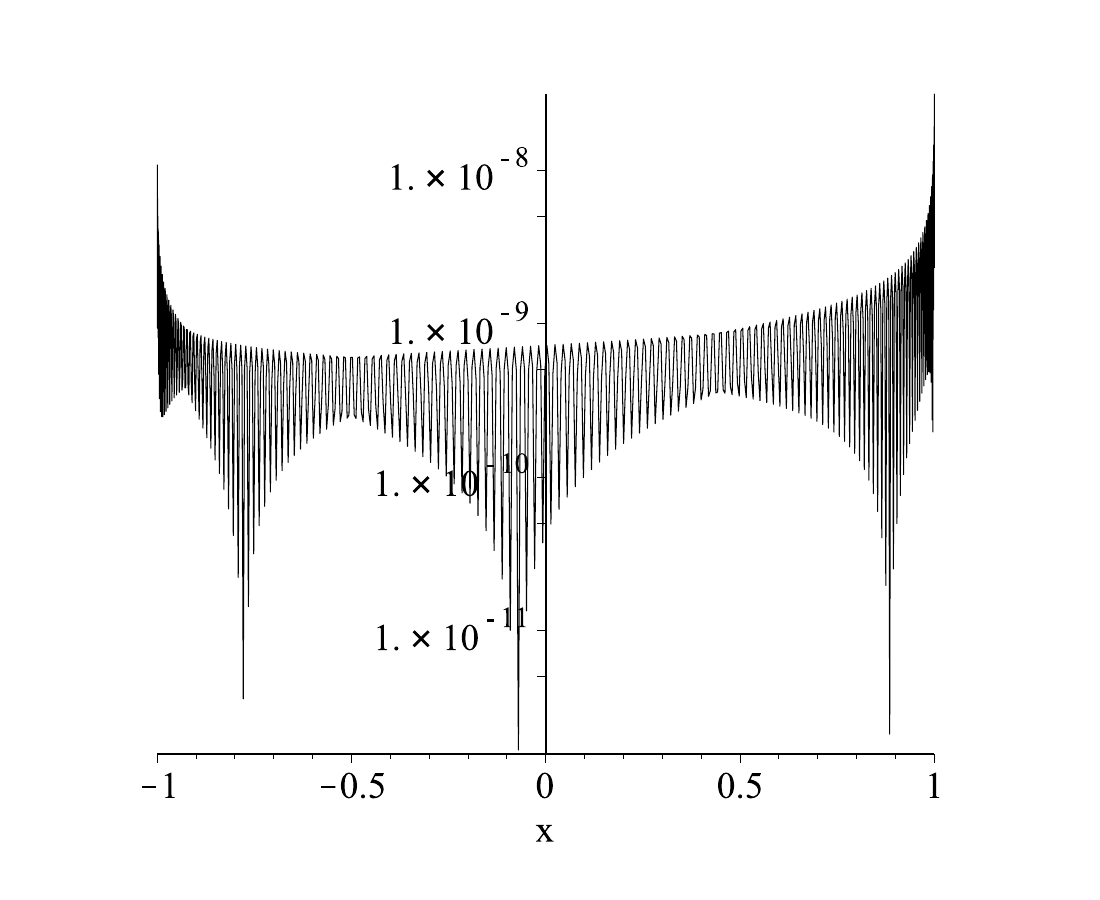}}
	\caption{Left : Green: error norm of the Chebyshev polynomial series, truncated at $n=N$, when the norm is the maximum over the \emph{middle} of the interval for  $u(x)=(1+x/2) (1 - x^{2})^{2} \, \ln(1 - x^{2})$. The green dashes are the reference line, $30/N^{5}$. The Chebyshev coefficients (not shown) also exhibit fifth order convergence, matching the power law exponent of the \emph{interior} error, for this $u(x)$. Black: same except the error norm is computed over the \emph{whole} of the interval, the usual norm. The black dashes are a graph of $1/N^{4}$.  Right : The errors versus $x$ of the Chebyshev polynomial series for $u(x)$, truncated after $N=50$ terms, that is, $E^{T}_{N}(x) \equiv |u(x) - u_{N}(x)|$. }
	\label{FigXXX}
\end{figure}

\begin{lemma}~\label{LemmaBD}
	For $\kappa \geq 2 $ and $\vartheta\in\mathbb{N}$, then
	\begin{equation}\label{EqBasz}
	 \sum_{n=N+1}^{\infty} \, \frac{ \ln^{\vartheta} (n) }{ n^{\kappa}  } \, \sim  \,\frac{\ln^{\vartheta}( N)}{(\kappa - 1) N^{\kappa-1}}.
	\end{equation}
\end{lemma}

The lemma is proved in \cite{BaszenskiDelvos88} when $\vartheta=0$ and in \cite{Zhang2021} when $\vartheta\in \mathbb{N}$.

Fig.~\ref{FigXXX} on the left side shows this steep rise in error near the endpoints by comparing two different norms. The upper solid (black) curve, falling $1/N$ slower than the coefficients, is the usual maximum pointwise error:
\begin{equation*}
E_{N}^{T}=\max_{x \in [-1,1]}|u(x) - u_{N}(x)|. 
\end{equation*}

The lower curve, which decreases as rapidly as the coefficients $a_n$, represents the maximum error over the interval's \emph{interior}, excluding the neighborhoods of both endpoints,


\begin{eqnarray*}
	E_{N}^{T,interior}=\max_{x \in [-0.8, 0.8]}|u(x) - u_{N}(x)| .
\end{eqnarray*}


Instead of plotting norms versus truncation as on the left side of Fig.~\ref{FigXXX},
a direct confirmation of the large errors in narrow boundary layers at the endpoints can be obtained by plotting errors versus $x$ as is shown on the right side of Fig.~\ref{FigXXX}.

	%

\begin{theorem}[Error in truncation of infinite series]
	Suppose we truncate each of the three series to a polynomial of degree $N$, as given in \eqref{eqn:trunCheb}, \eqref{eqn:trunDiff} and \eqref{eqn:trunQuad}.
	Then the error estimates of the truncated series in $L_{\infty}$ norm are presented in the following.
	\begin{enumerate}[(i).]
		\item For the Chebyshev series,
		\begin{eqnarray}
			E_{N} & \equiv & \max_{x \in [-1, 1]} | u(x) - u_{N}(x) | \sim \mathcal{O}\left(|A(N)|/N^{2 \varphi}\right), \qquad N \rightarrow \infty.
		\end{eqnarray}
		
		\item For the difference basis, 
		\begin{eqnarray}
			E^{diff}_{N} & \equiv & \max_{x \in [-1, 1]} | u(x) - u^{diff}_{N}(x) | \sim \mathcal{O}\left(|A(N)|/N^{2 \varphi}\right), \qquad N \rightarrow \infty .
		\end{eqnarray}
		
		\item For the quadratic-factor basis,
		 \begin{eqnarray}
			E^{quad}_{N} & \equiv & \max_{x \in [-1, 1]} | u(x) - u^{quad}_{N}(x) | \sim \mathcal{O}\left(|A(N)|/N^{2 \varphi-1}\right), \qquad N \rightarrow \infty. 
		\end{eqnarray}
The $A(N)$ is given in Theorem \ref{Thm1}.
	\end{enumerate}
\end{theorem}

\begin{proof}
: The error in the Chebyshev series follows from the discussion preceding the theorem. To prove the remaining propositions, note that the coefficients of the latter two expansions match up to degree $N-2$ when expanded as Chebyshev series. However, the difference relations in Theorem~\ref{Thabc} show that, with $b_{N-1}=b_{N}=c_{N-1}=c_{N}=0$, 
\begin{eqnarray}~\label{Eqtrunc3}
	u_{N}^{diff}(x) &= &\sum_{n=0}^{N-2} a_{n} T_{n}(x) + b_{N-3} T_{N-1}(x) \, + b_{N-2} \, T_{N}(x), \nonumber \\
             &= & u_{N}(x) + \left( b_{N-3} - a_{N-1} \right) T_{N-1}(x) \, +
	\left( b_{N-2} \, - a_{N} \right) T_{N}(x),   \nonumber\\
&= & u_{N}(x) +  b_{N-1} T_{N-1}(x) \, +
	\, b_{N} \, T_{N}(x) . 
\end{eqnarray}

Now we know from Theorem \ref{Thm1} that $b_{n} \sim \mathcal{O}\left(A(n)/n^{2\varphi}\right) $. This implies that $b_{N-1}$ and $b_N$ are proportional to the same power of $N$ as the error in the truncated Chebyshev series. It follows that the error in the truncated series on the difference basis has the same rate of convergence as the truncation of the Chebyshev series. 


To prove the final proposition, observe that the truncated series on the quadratic-factor basis can be written 
\begin{eqnarray}
	u^{quad}_{N}(x)&=&\sum_{n=0}^{N-2} c_{n} (1 - x^{2}) T_{n}(x)  \nonumber  \\
	& = & \left( \dfrac{1}{2} c_{0}  - \dfrac{1}{4} c_{2} \right)  T_{0}(x)  + 
	\left( \dfrac{1}{4} c_{1} - \dfrac{1}{4} c_{3} \right)T_{1}(x)
	+ 
	\left(- \dfrac{1}{2} c_{0} + \dfrac{1}{2} c_{2} - \dfrac{1}{4} c_{4} \right) T_{2}(x) \nonumber \\
	&+ & \sum_{n=3}^{N-2}    \left( \dfrac{1}{2} c_{n} - \dfrac{1}{4} c_{n-2} - \dfrac{1}{4} c_{n+2} \right) T_{n}(x)  
	-  \dfrac{1}{4} c_{N-3}  T_{N-1} (x) - \dfrac{1}{4} \, c_{N-2} T_N(x) \nonumber\\
	& = & \sum_{n=0}^{N} \, a_{n} T_{n}(x)    +\left(- \dfrac{1}{4} c_{N-2} - a_{N-1} \right) T_{N-1} (x) +\left(-\dfrac{1}{4} c_{N-2}-a_N   \right) T_N(x)  \nonumber  \\
	& = & u_N(x) -\dfrac{1}{2} c_{N-1} T_{N-1} (x) -    \dfrac{1}{2} c_{N} T_N(x).
\end{eqnarray}

Theorem \ref{Thabc} shows that $c_{N-1}$ and $c_{N}$ are $\mathcal{O}(A(N)/N^{2 \, \varphi-1})$. This is larger than the error in the truncated Chebyshev series by a factor of $N$. Thus this is the magnitude of the error in the truncated quadratic-factor basis.  
$\blacksquare$
\end{proof}

\begin{figure}[h]
	\centerline{\includegraphics[scale=0.5]{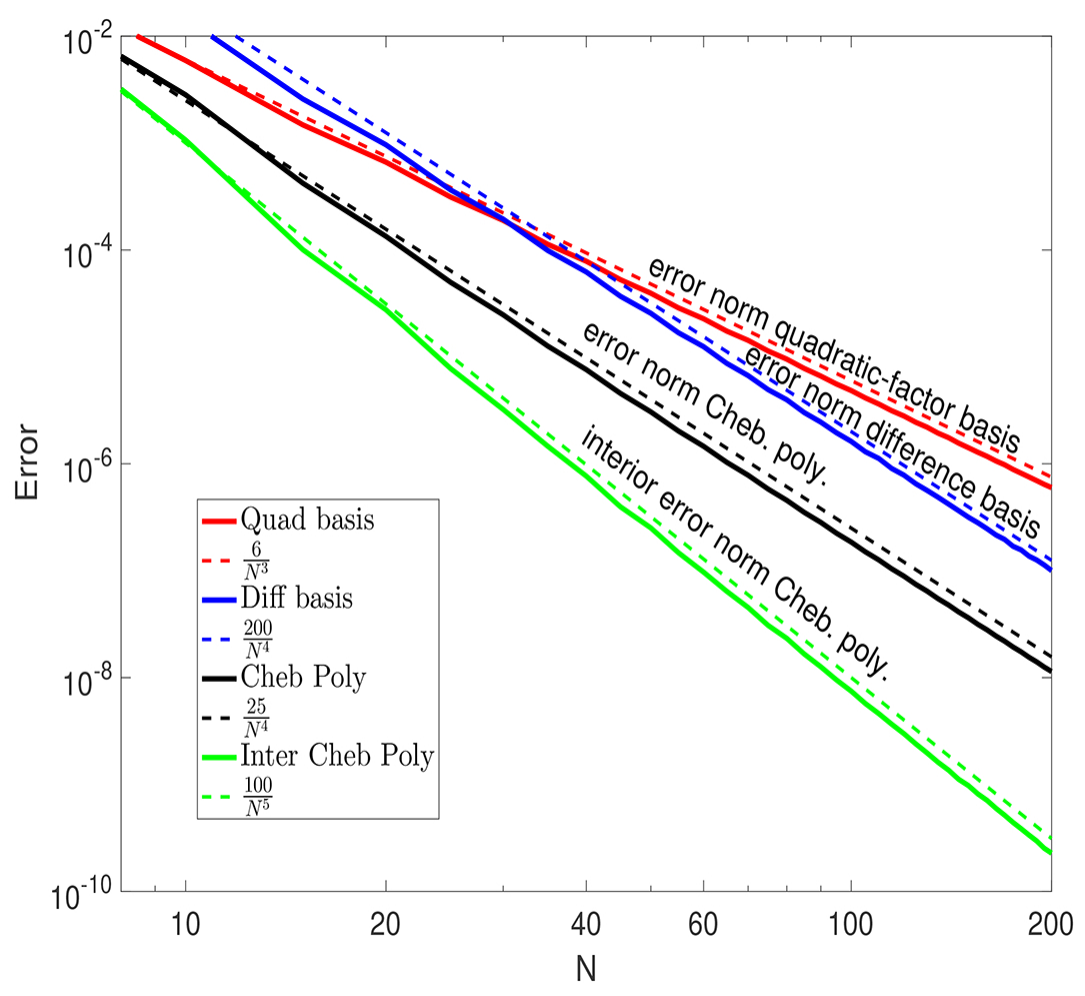}}
	\caption{Maximum pointwise errors ($L_{\infty}$ norm) in the truncated infinite series in three 
		different basis sets for various truncations $N$ for the typical example 
		$u(x)=(1+x/2)  (1- x^{2})^{2} \, \ln(1-x^{2})$. The upper solid line (red) is the error norm for 
		the quadratic-factor basis; the dashed red line is $6/N^{3}$. The solid blue curve is the 
		error norm for the difference basis; the dashed blue line is $200/N^{4}$.
		The solid black curve is the error norm for the truncation of the standard Chebyshev series; The black dashed curve is $25/N^{4}$. The green solid curve is the maximum pointwise error for $x \in [-1/2, 1/2]$, the interior of the interval $x \in [-1, 1]$; the dashed green curve is $100/N^{5}$.}
	\label{FigYYY}
\end{figure}

Fig.~\ref{FigYYY} confirms the expected rates of decay for an arbitrary but representative example.


\vspace{-0.4cm}
\section{Equivalence Theorem}

\begin{theorem}[Dirichlet-Enforcing Basis Equivalence]~\label{ThEquiv}
	If two polynomial approximations, constrained to satisfy homogeneous Dirichlet boundary conditions, are determined by the same set of interpolation constraints or least squares conditions, then the approximations are identical and must have identical errors, that is,
	\begin{eqnarray}~\label{Hump}
		u_{N}^{diff}(x) & = &  \,u^{quad}_{N}(x).
	\end{eqnarray}
	
\end{theorem}

\begin{proof} : By definition, $u^{diff}_{N}(x)$ is a polynomial 
of degree $N$ which is zero at both endpoints. The Fundamental Theorem of Algebra
asserts that any polynomial can be written in factored form. Therefore
\begin{eqnarray}
	u^{diff}_{N}(x) = (1 - x^{2}) \,p_{N-2}(x),\quad  N\geq 2,
\end{eqnarray}
where $p_{N-2}(x)$ is a polynomial of degree $(N-2)$. This is identical in form 
to $u^{quad}_{N}(x)$. If, for example, we determine the approximations by $N-1$ interpolation conditions, 
these constraints uniquely determine
$p_{N-2}(x)$ as the interpolant of $v(x)$, the same as for $u^{quad}_{N}(x)$. 
Therefore $u_{ N}^{diff}(x)  =   \,u^{quad}_{N}(x)$ for interpolation. The argument extends to 
any other reasonable mechanism to determine the approximations provided the same
conditions are applied to both $u^{diff}_{N}(x)$ and $u^{quad}_{N}(x)$. 
$\blacksquare$
\end{proof}

This equivalence theorem greatly simplifies error analysis. However, we  have already shown  
that the \emph{coefficients} $b_{n}$ and $c_{n}$ are \emph{different}. Furthermore, 
the error of an unconstrained series of Chebyshev polynomials is different 
from that of the constrained approximations.

\section{Interpolation \& Aliasing Errors in Chebyshev Polynomial Coefficients}

\subsection{Grids and uniqueness}~\label{SecGrid}

There are two canonical interpolation grids associated with Chebyshev polynomials. The ``roots" grid is
\begin{equation}
	x_k = - \cos\left(\frac{2k+1}{2N + 2}\pi \right)
	, k= 0,1,\,\ldots,\,N.  \qquad \mbox{[``Chebyshev-Gauss'' Grid]}
\end{equation}
The ``endpoints-and-extrema" or ``Lobatto" grid is 
\begin{equation}
	x_k = - \cos\left( \frac{k}{N} \pi \right)
	,  k= 0,1,\,\ldots,\,N.  \qquad \mbox{[``Chebyshev-Lobatto'' Grid]}
\end{equation}
If the Lobatto grid is chosen, then the interpolating polynomial must be 0 at $x=\pm 1$ in order to satisfy the interpolation condition at the endpoints. It follows that whether we represent the interpolated polynomial using Chebyshev polynomials, the difference basis, or the quadratic-factor basis, we always obtain the same polynomial. 

In contrast, if the interpolation points are those of the roots grid, which does not include the endpoints, then standard Chebyshev polynomial interpolation gives an interpolating polynomial which is not exactly equal to 0 at the endpoints. If we use either the quadratic-factor basis or the difference basis, the result, by the Polynomial Factorization Theorem, can be written in the form
\begin{eqnarray}
	u_{N}^{I,con}(x) = (1-x^{2}) v_{N-2}^{I,con}(x),
\end{eqnarray}
where $u_{N}^{I,con}(x), v_{N-2}^{I,con}(x)$ are Chevbyshev interpolants on Chebyshev-Lobatto grids for the function $u(x)$ and $v(x)=u(x)/(1-x^{2})$ respectively and they satisfy the homogeneous Dirichlet boundary conditions.
Thus, there are \emph{two} distinct interpolants on the \emph{roots} grid, these being the Chebyshev interpolant (lacking zeros at the endpoints) and the difference-and-quadratic-factor interpolant (which vanishes at both endpoints by construction). In contrast, the interpolant 
on the \emph{Lobatto} (endpoint-including) grid is always \emph{unique}.

\subsection{Aliasing errors in the Chebyshev coefficients of the interpolant}

The Chebyshev coefficients of both the interpolant, $a_{n}^{I}(N)$, and of the infinite series $a_{n}$ can be computed by Gauss-Chebyshev quadrature as given on pg. 99 of \cite{Boyd99z}. When the number of quadrature points $M$ is equal to $N+1$, then the coefficients are the result from interpolation; the coefficients of the infinite series are $a_{n}=\lim\limits_{M \rightarrow \infty, \mbox{fixed} N}a_{n}^{I}(N)$. But what is the relationship between series and interpolant coefficients for finite $N$? The following provides an answer.

\begin{theorem}[Aliasing formula for Chebyshev coefficients]~\label{Thalerr}
	Let $u(x)$ be Lispchitz continuous on $[-1,1]$ and let $u^{I}_{N}(x)$ be its Chebyshev interpolant 
	\begin{equation}
		u^{I}_{N}(x) = \frac{1}{2} a_{0}^{I} + \sum_{n=1}^{N} a_{n}^{I} T_{n}(x),
	\end{equation}
which is obtained by choosing the Chebyshev-Gauss grids as interpolation points.
Let $a_{n}$ (without superscript) denote the coefficients of the infinite series
	\begin{equation}\label{ChebExpan}
		u(x) = \frac{1}{2} a_{0} + \sum_{n=1}^{\infty} a_{n} T_{n}(x).
	\end{equation}
	Then, one has
	\begin{enumerate}[(i)]
		\item \begin{equation}
			\begin{aligned}
			a_{n}^{I} &= a_n + \mathcal{E}_{n}, \label{Oldfour5.17}\\
			\mathcal{E}_{n} & =    \sum_{j=1}^\infty \left( a_{-n+2j(N+1)}  +  a_{n+2j(N+1)}
			\right) (-1)^j, \;\,  n=0,1, \ldots, N.
			\end{aligned}
		\end{equation}	
		\item
		\begin{eqnarray}
			&&a_{n}^{I}  \approx  a_{n} - a_{2N+2-n} - a_{2N+2+n} + \mathcal{O}(a_{3N}), \qquad n=0, \ldots ,N ,\label{CoeffsAlias}\\
			&&a_{\frac{1}{2}(N+1)}^{I}  \approx  a_{\frac{1}{2}(N+1)} - a_{\frac{3}{2}(N+1)} + \mathcal{O}(a_{\frac{5}{2} (N+1) }), \nonumber \\
			&&a_{N+1-m}^{I}  \approx  a_{N+1-m} - a_{N+1+m}  + \mathcal{O}(a_{3N}), \, \, \, m\in \mathbb{N}_{+}, \; m \ll N\nonumber.
			\end{eqnarray}
	\end{enumerate}
	
\end{theorem}

\begin{proof} : 
The first proposition was proved by Fox and Parker \cite{FoxParker68}. 
The second comes from specializing to particular ranges in degree and then making obvious approximations. 
$\blacksquare$
\end{proof}


\begin{theorem}~\label{Th6}
	Suppose that the Chebyshev coefficients in \eqref{ChebExpan} for large  $n$ are 
	\begin{eqnarray*}
		a_{n} \sim A\frac{ \ln^{\vartheta}(n) }{ n^{\kappa} },\quad n\gg 1,\quad \kappa>0, \quad  \text{A is a constant and}\; \vartheta\in\mathbb{N} .
	\end{eqnarray*}
Then one has the following estimates :
	\begin{enumerate}[(i).]
		\item For small degree $n$, the aliasing error in Chebyshev coefficients is
		\begin{equation}~\label{Eq45}
\mathcal{E}_{n} \sim  -\frac{ A  \ln^{\vartheta}(2N) } {2^{\kappa-1} N^{\kappa}},	
		\end{equation}
	and the relative error is
	\begin{eqnarray}
		\frac{ | \mathcal{E}_{n}| }{|a_{n}|}  &  \sim &   \frac{1}{2^{\kappa-1} }  \frac{  n^{\kappa} }{  N^{\kappa}  } \frac{ \ln^{ \vartheta } (2N) } { \ln^{\vartheta} (n)}.
	\end{eqnarray}
	
		Specially, if the coefficients $a_{n}$ are well-approximated by the power law $A/n^{\kappa}$, for small degree $n$ such that $ n \ll   \frac{2}{\kappa} \, (N+1)$,
	then
		\begin{eqnarray}
			\mathcal {E}_{n} &\sim & \frac{A}{2^{\kappa-1}} \frac {1}{ N^{\kappa}}\sum_{j=1}^{\infty}. \frac{(-1)^{j}}{j^{\kappa}},\qquad
			\frac{ | \mathcal{E}_{n}| }{|a_{n}|}    \leq     \frac{1}{2^{\kappa-1} }  \frac{  n^{\kappa} }{  N^{\kappa}  } .
		\end{eqnarray}
		
		\item \label{AliasingSmallm}For $n=N+1-m$ when $m$ is a small positive integer,
		the relative error is
		\begin{eqnarray}
			\frac{ | \mathcal{E}_{N+1-m}| }{|a_{N+1-m}|} \sim 1 + \mathcal{O}\left( \frac{ \kappa m}{N+1} \right).
		\end{eqnarray}
	\end{enumerate}
\end{theorem}
\begin{proof}
: Substituting the coefficients into the terms in the error sum gives
\begin{eqnarray*}
&&	a_{2j(N+1) \pm n}    \sim   \frac{A \ln^{\vartheta} (2j N) }{ (2 j (N+1) \pm n)^{\kappa}}
	 \sim \frac{A\ln^{\vartheta}(2jN)}{2^{\kappa} (N+1)^{\kappa}}    \frac{1}{ j^{\kappa} \,  \left\{ 1 \pm n/(2 j (N+1))  \right\}^{\kappa} }, \;\\
&&	a_{2j(N+1)-n} +  a_{2j(N+1)+n}
	\approx  \frac{A  \ln^{\vartheta} (2jN)  }{2^{k-1} (N+1)^{\kappa}}    \frac{1}{ j^{\kappa} \,   }  \, \left\{ 1 + \mathcal{O}\left(\frac{\kappa^{2} n^{2}}{4 j^{2} (N+1)^{2}}\right) 
	\,\right\}.
\end{eqnarray*}
The asymptotic expression (\ref{Eq45}) then follows.

%

If we assume $n$ is sufficiently large that $a_{n}\sim A/ n^{\kappa}$, then the relative coefficient error follows 
immediately upon invoking
\begin{eqnarray}
	\sum_{j=1}^{\infty} \frac{ (-1)^{j+1} }{j^{\kappa}} < 1, \qquad \forall \kappa > 0.
\end{eqnarray}

To prove the second proposition, substitute the coefficients decay law into the equation \eqref{CoeffsAlias} from 
the general aliasing theorem, Theorem~\ref{Thalerr}, using asymptotic tools, the item \eqref{AliasingSmallm} will come by.
$\blacksquare$
\end{proof}

The second proposition implies that coefficients whose degree is near the aliasing limit, $n=N$, are badly in error. When the Chebyshev coefficients decay slowly, as $\mathcal{O}(\ln^{\vartheta}(n)/n^{\kappa})$, 
$a_{N+1-m} \approx a_{N+1+m}$; the relationship $a_{N+1-m}^{I}  \approx  a_{N+1-m} - a_{N+1+m} $ 
implies strong cancellation so that
\begin{eqnarray}
	| a_{N+1-m}^{I} | \ll |a_{N+1-m} |.
\end{eqnarray}
and the relative error is near 100 \%. When $\vartheta=0$, a log-log plot of the interpolation points $|a_{n}^{I}|$ is a straight line for intermediate $n$, but dives to small values as $n \to N$, curving downward below the line. For general $\vartheta \in\mathbb{N}_{+}$, the plot of interpolation coefficients $|a_{n}^{I}|$ still curl downward much more than the curve of Chebyshev coefficients $|a_{n}|$ as $n\to N$. 

The first proposition shows that in contrast, low degree coefficients can be computed with 
small relative error, but for fixed degree $n$, the relative error falls with $N$ as $\ln^{\vartheta}(N)/N^{k}$, the same decay rate as the coefficients. In words, if $a_{n}$ diminishes with $\ln^{\vartheta}(n)/n^{k}$ the coefficient can be computed as the corresponding coefficient of the $(N+1)$-point interpolant with a relative error that is order $k$ in $N$ by a factor of $\ln^{\vartheta}(N)$.


\section{Interpolants and Interpolation Errors with Dirichlet Boundary Conditions}

\subsection{Interpolants and their similarities and differences}

Because the Lobatto grid includes the endpoints, the standard, unconstrained Chebyshev interpolant is zero at both endpoints for any function satisfying $u(\pm 1)=0$. As noted in Sec.~\ref{SecGrid}, the interpolant on the Lobatto grid is unique and therefore:
\begin{eqnarray}
	u_{N}^{Cheb,Lob,I}(x) = u_{N}^{diff,Lob,I}(x) = u_{N}^{quad,Lob,I}(x).
\end{eqnarray}

So, let us turn to the roots grid. Define $v(x) \equiv u(x)/ (1-x^{2})$ as before. There exists a polynomial of degree $(N-2)$, which we will denote by  $v_{N-2}^{Cheb,I}(x)$, that interpolates $v(x)$ at all of the points on the $(N-1)$-point roots grid.

\begin{theorem}[Interpolants on the Roots Grid]
	Suppose that $u(x)$ satisfies Dirichlet boundary conditions $u(\pm1)=0$
and the $(N-1)$-point Chebyshev interpolant of $v(x)$ is
	\begin{eqnarray}
		v^{Cheb,I}_{N-2}(x) & =& \sum_{n=0}^{N-2}   \tilde{c}_{n}^{I} \, T_{n}(x).
	\end{eqnarray}
	Compute $u^{quad,I}_{N}(x)$ by $(N-1)$-point interpolation of $u(x)$ where
	\begin{eqnarray}
		u^{quad,I}_{N}(x) & =& \sum_{n=0}^{N-2}   c_{n}^{I} \, (1-x^{2}) \, T_{n}(x).
	\end{eqnarray}
	Similarly compute $u^{diff,I}_{N}(x)$ by $(N-1)$-point interpolation where
	\begin{eqnarray}
		u_{N}^{diff,I}(x) = \sum_{n=0}^{N-2} b_{n}^{I} \left\{ T_{n+2}(x) \, - \, T_{n}(x) \right\} .
	\end{eqnarray}
	Then, leads to
	\begin{eqnarray}
		u_{N}^{quad,I}(x) & = & (1- x^{2}) \,v^{Cheb,I}_{N-2}(x) ,\\
		u_{N}^{diff,I}(x) & = & (1- x^{2}) \,v^{Cheb,I}_{N-2}(x), \label{eq118} \\
		u_{N}^{diff,I}(x) & = & u_{N}^{quad,I}(x), \label{eq119} \\
		c_{n}^{I}  =  \tilde{c}_{n}^{I} , \;i&= &0,1, \, \ldots,\, N-2,
	\end{eqnarray}
and 
\begin{eqnarray}
	b_{0}^{I} & = & \dfrac{c_{2}^{I} - 2 c_{0}^{I}}{4} ,\;
	b^{I}_{n}  =  \dfrac{c_{n+2}^{I} -  c_{n}^{I}}{4}, \, \, \,  n=1,2, \ldots ,(N-4), \;
	b^{I}_{N-3}  =  - \dfrac{c_{N-3}^{I}}{4} ,\;
	b^{I}_{N-2}  =  - \, \dfrac{   c_{N-2}^{I}}{4}.
\end{eqnarray}
	
\end{theorem}
\begin{proof}
: The interpolation conditions for $u(x)$ in the quadratic-factored basis are
\begin{eqnarray}~\label{Burst}
	u(x_{j}) &  =  &  u^{quad, I}_{N}(x_{j})
	 = \sum_{n=0}^{N-2}   c_{n}^{I} \, (1-x_{j}^{2}) \, T_{n}(x_{j}).
\end{eqnarray}
The same for $v(x)$ multiplied by $(1-x^{2})$ are
\begin{eqnarray}~\label{Bomb}
	(1- x_{j}^{2})  v(x_{j}) &  =  &  (1 - x_{j}^{2}) \,v^{Cheb,I}_{N-2}(x_{j})
	 = \sum_{n=0}^{N-2}  \tilde{c}_{n}^{I} \, (1-x_{j}^{2}) \, T_{n}(x_{j}).
\end{eqnarray}
The left-hand side of (\ref{Bomb}) is $u(x_{j})$. The right-hand side  is identical in 
form to the interpolant of $u(x)$ by $u_{N}^{quad,I}(x)$. Therefore $\tilde{c}_{n}^{I}=c_{n}^{I}$
from which follows $u_{N}^{quad,I}(x) =  (1- x^{2}) \,v_{N-2}^{Cheb,I}(x)$.
The second and third lines, (\ref{eq118}) and (\ref{eq119}), follow from the Equivalence Theorem~\ref{ThEquiv}.
The formulas for the $b_{n}^{I}$ follow from the difference equations in Proposition 6 of Theorem~\ref{Thabc}.
$\blacksquare$
\end{proof}


\subsection{Interpolation Errors and Error Norms}

Suppose that the Chebyshev polynomial coefficients $a_{n}$ of a function $u(x)$ are decreasing as 
\begin{eqnarray}
	a_{n} \sim A \ln^{\vartheta}(n) / n^{\kappa},
\end{eqnarray} where here $\kappa=2 \varphi+1>0, \vartheta\in \mathbb{N}$. 
The error in the Chebyshev interpolant of $u(x)$ is expected to be $\mathcal{O}(\ln^{\vartheta}(N)/N^{\kappa})$ on the interior of the interval, slowing to $\mathcal{O}( \ln^{\vartheta}(N) /N^{\kappa-1})$ in the endpoint boundary layers.

The Chebyshev polynomial coefficients of $v(x) \equiv u(x)/(1-x^{2})$ converge more slowly than those of $u(x)$ by a factor of about $n^{2}$ (Tuan and Elliott \cite{TuanElliott72}). Define 
\begin{eqnarray}
	E_{N}^{v}(x) \equiv |v(x) - v_{N-1}^{I,con}(x)| .
\end{eqnarray}
It follows that $E_{N}^{v}(x)$ will be $\mathcal{O}(\ln^{
\vartheta}(N)/N^{\kappa-1})$ on the interior of the interval. To obtain the corresponding error in $u(x)$, we must multiply by the factor of $(1-x^{2})$ which is the ratio of $u(x)$ to $v(x)$, that is 
\begin{eqnarray}
	E_{N}^{u} (x) = (1-x^{2}) E_{N}^{v}(x).
\end{eqnarray}
It follows that 
\begin{eqnarray}
	\max_{x \in [-1,1]} \left (E_{N}^{u}(x) \right ) \sim \frac{A\ln^{\vartheta}(N)}{N^{\kappa-1}}.
\end{eqnarray}

\begin{figure}
	\centerline{\includegraphics[scale=0.6]{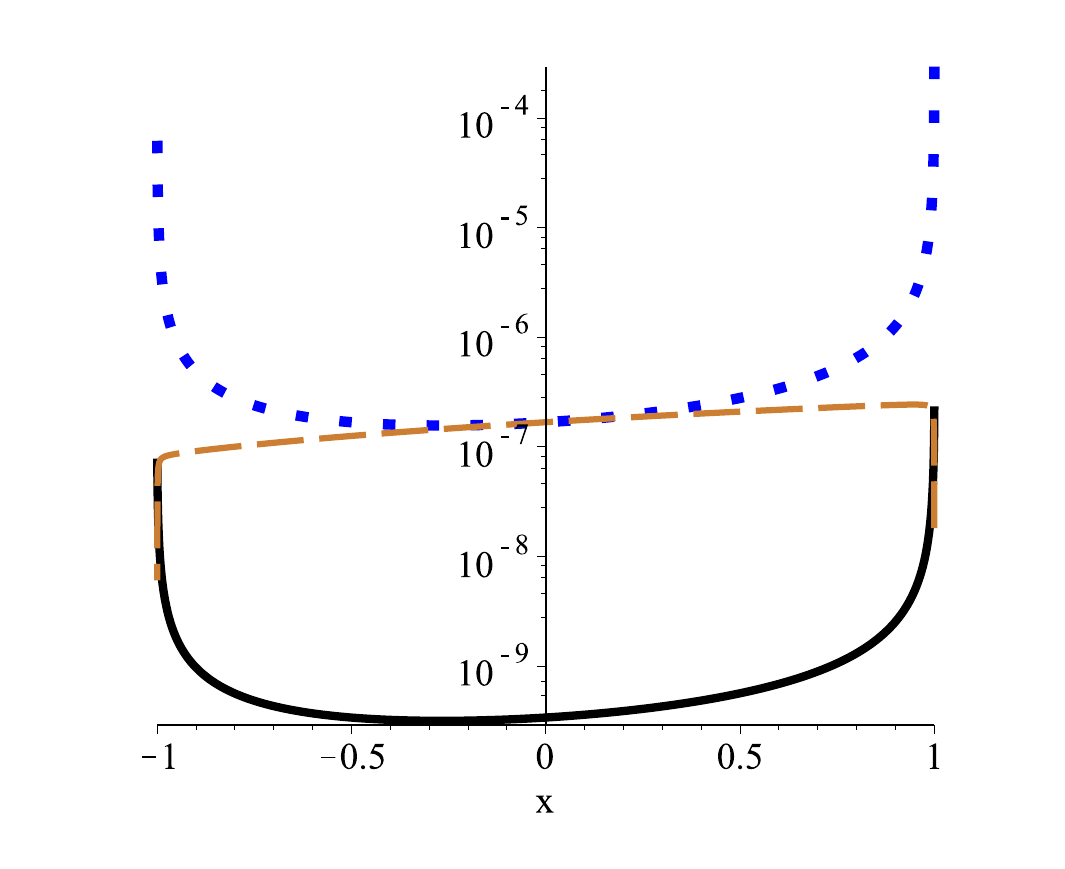}}
	\caption{Errors versus $x$ for interpolation of 
		$u(x)=(1+x/2)  (1- x^{2})^{2} \, \ln(1-x^{2})$, the same function as employed in the previous figure, by means of 100 interpolation points on the roots grid.  Top (blue dots): $|v(x) - v_{N}^{Cheb,I}(x)|$, the interpolant of $v(x)=u(x)/(1-x^{2})$.   Bottom (solid black curve): $|u(x)  - u_{N}^{Cheb,I}(x)|$, the error in the classic Chebyshev interpolation on the roots grid.  The dashed gold line is the error for
		the quadratic-factor basis, the error in $u_{N}^{quad,I}(x)= (1-x^{2}) v_{N-2}^{Cheb,I}(x)$; this is 
		identical with the error in the  difference basis since (for interpolation) $u_{N}^{quad,I}(x)=u_{N}^{diff,I}(x)$.
	}	\label{FigOP279_three_interp_err_vs_x}
\end{figure}

As explained in Chapter 2 of \cite{Boyd99z}, the error in truncating the infinite Chebyshev series by discarding all terms of degree $(N+1)$ and higher can be bounded rigorously by the sum of the absolute values of the neglected coefficients:
\begin{eqnarray}
	E_{N}(x) \equiv |u(x) - u_{N}(x)| \leq \sum_{n=N+1}^{\infty} |a_{n}| .
\end{eqnarray}
Chebyshev interpolation on either the roots or Lobatto grids is bounded by \emph{twice} the sum of the absolute 
values of the neglected coefficients:
\begin{eqnarray}
	E^{I}_{N}(x) \equiv |u(x) - u^{I}_{N}(x)| \leq \, 2 \, \sum_{n=N+1}^{\infty} |a_{n}| .
\end{eqnarray}
It is difficult to make more precise statements; for $u(x)= T_{\frac{3}{2}(N+1)}(x)$, for example,
\begin{eqnarray*}
& 	u_{N+1}(x)  = 0, \quad & E_{N+1}(x) = T_{\frac{3}{2}(N+1)}(x), \; \\
& 	u^{I}_{N}  =  - T_{\frac{1}{2}(N+1)}(x), \quad & E^{I}_{N+1}(x) = T_{\frac{3}{2}(N+1)}(x) + T_{\frac{1}{2}(N+1)}(x).
\end{eqnarray*}

Nevertheless, it follows that $E^{I}_{N}(x)$ is roughly double the error in truncating the infinite 
Chebyshev series and therefore its $L_{\infty}$ \emph{norm} is $\mathcal{O}(\ln^{\vartheta}(N)/N^{\kappa-1})$.

Because of the endpoint singularities, the usual nearly-uniform error for truncated Chebyshev series (or Chebyshev interpolants) of \emph{smooth} functions, analytic on the entire interval, is replaced by an error which is huge in boundary layers near each endpoint and smaller outside of these boundary layers by a factor of $\mathcal{O}(1/N)$ (bottom curve in 
Fig.~\ref{FigOP279_three_interp_err_vs_x}).

Applying this same reasoning to $v(x) \equiv u(x)/(1-x^{2})$ gives an error for $v(x)$ which is $\mathcal{O}(N)$ times as large as the error for $u(x)$ (Note that the order $\kappa$ of singularities for $v(x)$ is one less than for $u(x)$ and each decrease in $\kappa$ by \emph{one} reduces the order of the Chebyshev coefficients by \emph{two}). To get the approximation in the quadratic-factored basis for $u(x)$, we must multiply the Chebyshev series for $v(x)$ by $(1-x^{2})$. Similarly, the highly nonuniform error in $v(x)$ (blue dotted curve in Fig.~\ref{FigOP279_three_interp_err_vs_x}) musts be replaced, to get the error for interpolation of $u(x)$ by either of the constrained basis sets, by $(1-x^{2})$ times the error for $v(x)$. The zeros at the endpoints wipe out the boundary layers of large error in $v(x)$ to yield an error which is nearly uniform over $x \in [-1, 1]$ as shown by the gold dashed curve in Fig.~\ref{FigOP279_three_interp_err_vs_x}.

Plain classical Chebyshev interpolation, although no better than the other two basis sets 
in the $L_{\infty}$ norm, is superior because the pointwise Chebyshev interpolation is as 
bad-as-the-norm only in boundary layers  whereas the quadratic-factor and difference errors
are as bad as the norm over the entire interval.

\begin{figure}
	\centerline{\includegraphics[scale=0.6]{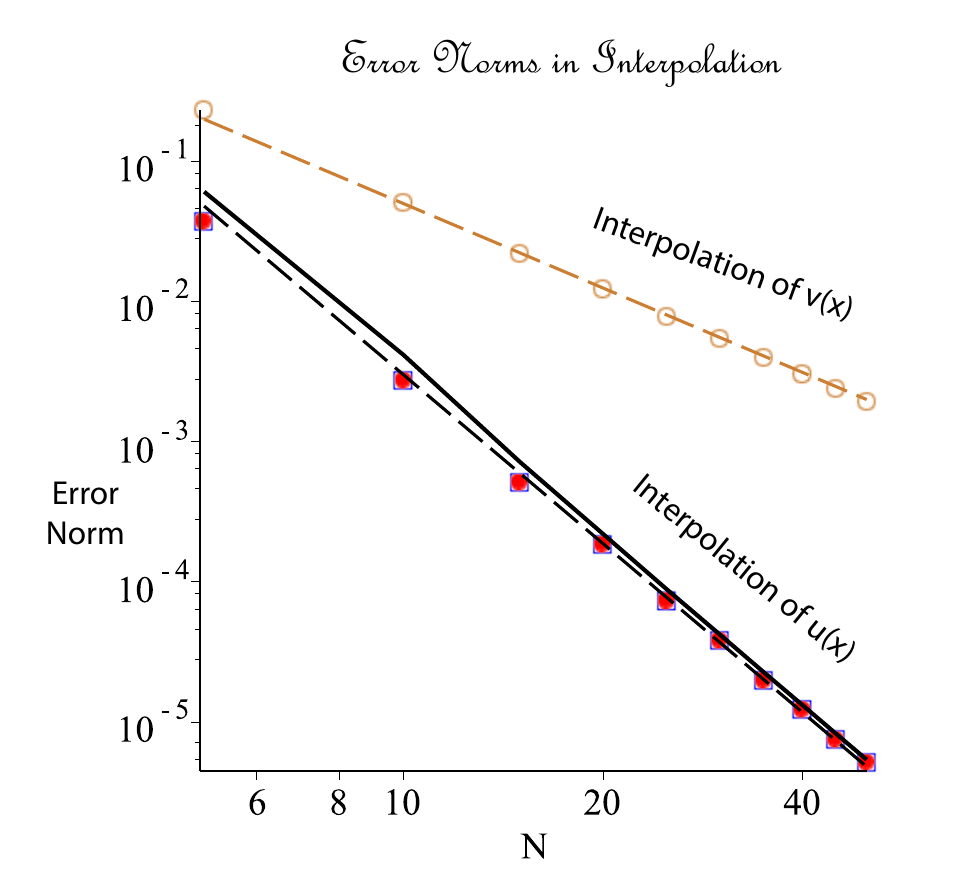}}
	\caption{Log-log plot of the error norms  versus the number of interpolation points $N$ for interpolation of $u(x)=(1+x/2)  (1- x^{2})^{2} \, \ln(1-x^{2})$, the same function as employed in the previous figure.  Top circle gold curve: error norm $\max_{x \in [-1, 1]}|v(x) - v_{N}^{Cheb,I}(x)|$, the interpolant of $v(x)=u(x)/(1-x^{2})$. The gold dashed curve is $5/N^{2}$, proportional to $1/N^{\kappa-3}$.  The bottom three curves, almost superimposed and hard to distinguish, are the error norms for the approximation of $u_{N}^{Cheb,I}(x)$ by the Chebyshev interpolant (black curve),  $u_{N}^{quad,I}(x)$ (blue boxes) and $u_{N}^{diff,I}(x)$ (red solid disks). The red disks are at the center of the blue boxes because $u_{N}^{quad,I}(x)$ is identically equal to $u_{N}^{diff,I}(x)$, as asserted by Theorem \ref{ThEquiv}. The black dashed curve is $30/N^{4}$, proportional to $1/N^{\kappa-1}=1/N^{2 \varphi}$ where $\kappa=5$. }
	\label{FigOP279_error_norm_diff_quad_Cheb_INTERP_MAPLE}
\end{figure}

Fig.~\ref{FigOP279_error_norm_diff_quad_Cheb_INTERP_MAPLE} displays error norms instead of 
pointwise errors. The close agreement between the dashed curves and the matching solid curves
confirms the theoretical predictions given above.

\subsection{Coefficients of Interpolants}

\begin{figure}
	\centerline{\includegraphics[scale=0.5]{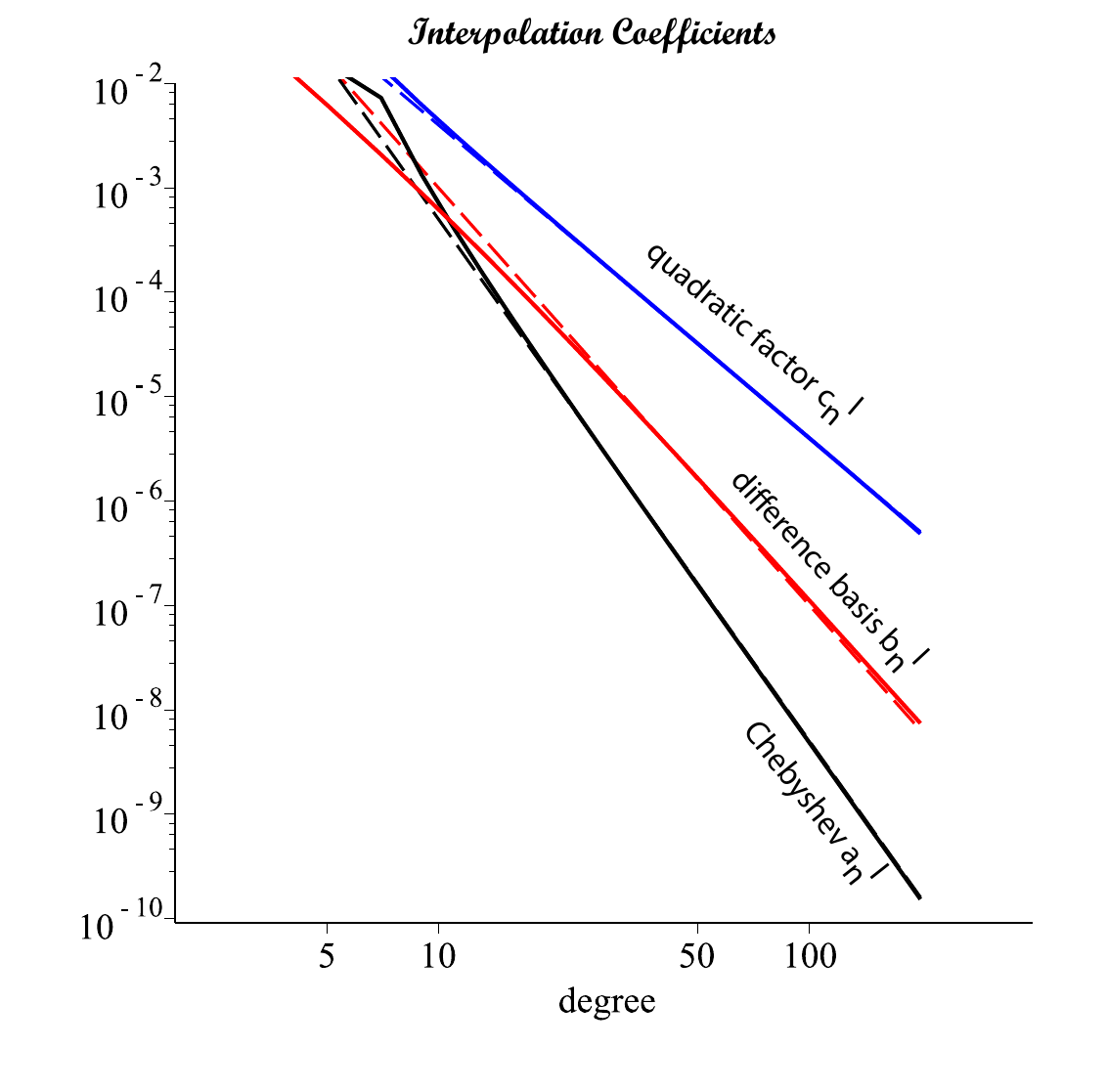}}
	\caption{Odd degree coefficients versus degree for interpolation of 
		$u(x)=(1+x/2)  (1- x^{2})^{2} \, \ln(1-x^{2})$, the same function as employed in the previous figure, by means of 400 interpolation points on the roots grid.  To minimize aliasing contamination, only the coefficients  up to degree 199 are plotted. Top (blue solid line): $c_{n}^{I}$, which are simultaneously the Chebyshev coefficients of $v(x)$ and also the coefficients of  $u(x)$ in the quadratic factor basis where $v(x)=u(x)/(1-x^{2})$. Middle (red solid curve): coefficients $b_{n}^{I}$ of the difference basis for $u(x)$. Bottom (solid black curve): the Chebyshev coefficients $a_{N}^{I}$ of the interpolant on the roots grid.  The dashed lines are proportional to $1/n^{3}$ (top/blue), $1/n^{4}$ (middle red) and $1/n^{5}$ (bottom/black). the parameter $\varphi=2$ while $\kappa=2 \varphi + 1 =5$.}
	\label{FigOP279_COEFFS_diff_quad_Cheb_INTERP_MAPLE}
\end{figure}

Fig.~\ref{FigOP279_COEFFS_diff_quad_Cheb_INTERP_MAPLE} shows how the coefficients vary. Even though the \emph{errors} of the difference basis and quadradic factor basis interpolants are 
\emph{identical}, their \emph{coefficients} obey \emph{different power laws}. The quadratic factor coefficients $c_{n}$ decay more slowly by one order than the difference basis coefficients $b_{n}$. The power laws for all three basis sets are the same as for truncation of the infinite series, so no further discussion will be given.

\vspace{-0.5cm}

\section{Least Squares}
Least squares is a third strategy that provides an alternative to interpolation and truncation of infinite series. Is it better? Worse? One complication is that least squares is actually a family of methods because the approximation varies with the choice of the inner product.

The next two subsections describe the basic methods with and without Lagrange multipliers. 
In the rest of the section, we shall analyze the least squares for three bases in turn. When the 
inner product is integration over the interval, we shall show that least squares yields approximations \emph{different} from interpolation and truncation for the constrained-to-vanish-at-the-endpoints basis sets.

\subsection{Least Squares without Lagrange Multiplier: General Basis}\label{Sec7_1}

The goal of least squares is to minimize the ``cost function"
\begin{eqnarray}
	\mathfrak{J} = \frac{1}{2} < u(x) - u_{N}(x) , u(x) - u_{N}(x) >,
\end{eqnarray}
where 
\begin{eqnarray}
	u_{N}(x)= \sum_{n=0}^{N} d_{n} \phi_{n}(x).
\end{eqnarray}
For the moment, the choices of inner product $< f(x), g(x)>$ and the basis functions $\phi_{n}(x)$ are unspecified.

\begin{proposition}
	Suppose that the cost function and $u_{N}(x)$ are as above. Define a $(N+1) \times (N+1)$ matrix $\mathbf{G}$ as the matrix with the elements
	\begin{eqnarray}  G_{mn} = < \phi_{m}(x), \phi_{n}(x) >, \qquad m,\, n=0,1, \ldots ,N.
	\end{eqnarray} 
	Define $\mathbf{f}$ as the vector with elements 
	\begin{eqnarray}
	  f_{n} \equiv < \phi_{n}(x), u(x) >,\quad n=0,1,\ldots, N.
	\end{eqnarray}
	Then $u_{N}(x)$ is the unique minimizer if and only if the spectral coefficients $d_{n}$ are the elements of 
	the  $(N+1)$-dimensional vector $\mathbf{d}$ which solves
	\begin{eqnarray}~\label{EqGsolve}
		\mathbf{G} \mathbf{d}=\mathbf{f}.
	\end{eqnarray}
\end{proposition}
\begin{proof} :
Substitute the series into the cost function and apply the condition for a minimum  that the derivatives of the cost function with respect to the coefficients $d_{j}$ are all zero. This gives
\begin{eqnarray}
	\frac{\partial \mathfrak{J}} {\partial d_{m}} = 0 =  - \, <u(x), \phi_{m}(x)> + \sum_{n=0}^{N} d_{n} < \phi_{m}(x), \phi_{n}(x) >,
\end{eqnarray}
which is the linear algebra problem (\ref{EqGsolve}). 
$\blacksquare$
\end{proof}

Let us suppose that the inner product is approximated by Gaussian quadrature with $N_{col}$ points. For the Chebyshev weight,
\begin{equation}
	\begin{aligned}
	< f(x), g(x)> \equiv \int_{-1}^{1} &\frac{1}{\sqrt{1 - x^{2}} } f(x) g(x) dx 
	 \approx   <f, g>_{Gq} \equiv  \frac{\pi}{N_{col}}  \sum_{j=1}^{N_{col}} \, f(x_{j}) g(x_{j}), \nonumber \\
\quad	&x_{j}  =  \cos\left( \, \frac{2 j - 1}{2 N_{col} }  \,\pi \, \right), \qquad j=1, 2, \ldots N_{col}. \nonumber
	\end{aligned}
\end{equation}
The quadrature approximation has all the properties to be an inner product, so we use $< \cdotp,\cdotp>_{Gq}$ as the inner product in the rest of the section. This inner product varies from interpolation (when $N_{col}=N+1$ as explained below) to integration over the interval in the limit $N_{col} \rightarrow \infty$.

Define a $N_{col} \times (N+1)$ matrix $\mathbf{H}$ whose elements are
\begin{eqnarray}
	H_{jn} \equiv \phi_{n}(x_{j}), \qquad j=1,2,\ldots, N_{col},  \, n=0,1,\ldots, N.
\end{eqnarray}
and let $\mathbf{u}$ denote the vector whose elements are the samples of $u(x)$,
\begin{eqnarray}
	u_{j} \equiv u(x_{j}),\, j=1,2,\ldots,N_{col}.
\end{eqnarray}
The interpolation problem is
\begin{eqnarray}~\label{EqLsinterp}
	\mathbf{H} \mathbf{d}^{I} = \mathbf{u}.
\end{eqnarray}
Here we have added a superscript to the vector of spectral coefficients because the solution to the interpolation problem is not necessarily the same as the solution that minimizes $\mathfrak{J}$. Note that the matrix problem is an overdetermined system, but still well-posed if $N_{col} > N+1$.

\begin{proposition} 
	The solution $\mathbf{d}$ of the least squares with an inner product using $N_{col}=N+1$ quadrature points 
	is identical to the solution $\mathbf{d}^{I}$  to the $N+1$ point interpolation problem.
	The matrices for least squares and interpolation are connected by
	\begin{eqnarray}
		\mathbf{G} = \frac{\pi}{N_{col}}\mathbf{H}^{T} \mathbf{H} ; \qquad  \mathbf{f} = \frac{\pi}{N_{col}}  \mathbf{H}^{T} \mathbf{u}.
	\end{eqnarray}
\end{proposition}

To prove this theorem, one can refer to the proof of Theorem 16 in \cite{BoydBook3}.

When the basis functions are orthogonal, the $n$-th element of the solution is independent of $N$ so long as $N \geq n$. The quadratic-factor and difference basis are not orthogonal, and the solution elements depend on $N$.

\subsection{Lagrange Multiplier Theorem: Equality of Minimizers }

When a ``cost" function $\mathcal{J}$ is to be minimized subject to the constraints $\Psi = 0$ and $\Omega=0$, it is very convenient to convert the problem to the unconstrained minimization of the modified ``cost" function 
\begin{eqnarray}
	J = \mathcal{J}\,  +\,  \lambda \, \Psi \, + \, \mu \, \Omega,
\end{eqnarray}
where $\lambda$ and $\mu$ are additional unknowns called  ``Lagrange multipliers" and $\Psi, \Omega$ denotes the boundary constraints $u(1)=0, u(-1)=0$ respectively. The conditions for a minimum, taking the original unknowns to be $\{ d_{j} \}$, are
\begin{eqnarray}
	\frac{\partial J}{\partial d_{j}} &=& 0 , \qquad j=0,1, \ldots, N, \\
	\frac{\partial J}{\partial \lambda} &=& \Psi = 0, \\
	\frac{\partial J}{\partial \mu} &=& \Omega = 0.
\end{eqnarray}
Here the constraints are $u(\pm1)=0$; expressed in terms of Chebyshev coefficients these are:
\begin{eqnarray}
	\sum_{j=0}^{N}  d_{j} & = &  0    \quad (\Leftrightarrow u(1)=0) ;\qquad 
	\sum_{j=0}^{N} \, (-1)^{j} \,   d_{j}  =  0  \quad (\Leftrightarrow u(-1)=0).
\end{eqnarray}

\begin{theorem}
	Consider two minimization problems. 
	\begin{enumerate}[(i)]
	\item . Suppose that $u_{N}^{con}(x)$ is a solution to the cost function 
	\begin{eqnarray}
		\mathcal{J} \equiv < u(x) - u_{N}^{con}(x), u(x) - u_{N}^{con}(x)>,
	\end{eqnarray}
where $u_{N}^{con}(x)$ is a polynomial of degree $N$ constructed so that $\Psi=0$ and $\Omega=0$, is satisfied independent of the remaining unknowns. For example, 
	\begin{eqnarray}
		u_{N}^{con}(x)= \sum_{n=0}^{N-2} b_{n} \left\{ T_{n+2}(x) - T_{n}(x) \right\}.
	\end{eqnarray}

	\item . Suppose that $u_{N}(x)$ is a solution to the cost function 
	\begin{eqnarray}
		J \equiv < u(x) - u_{N}(x), u(x) - u_{N}(x)> + \lambda \Psi + \mu \, \Omega, 
	\end{eqnarray}
where $u_{N}(x)$ is a polynomial of degree $N$, to be unconstrained-at-the-endpoints minimizer of the cost function. Then the two solutions $u_{N}^{con}(x)$ and $u_{N}(x)$ to the minimization problems $\mathcal{J}$ and $J$ respectively are identical. 
	
	
\end{enumerate}
\end{theorem}

\begin{proof}:
 Now the solution to the second minimization problem is forced to satisfy the constraint as well. At the minimum, $\Psi =\Omega = 0$, so the cost function reduces to 
\begin{eqnarray}
	J \equiv < u(x) - u_{N}(x), u(x) - u_{N}(x)>.
\end{eqnarray}
It follows that $u_{N}(x)$ and $u_{N}^{con}(x)$ both minimize $< u(x) - u_{N}^{approx}(x), u(x) - u_{N}^{approx}(x)>$ where $u_{N}^{approx}(x)$ is either $u_{N}(x)$ or $u_{N}^{con}(x)$. Therefore, $u_{N}(x) \neq u_{N}^{con}(x)$ if and only if $u_{N}^{approx}(x)$ is not unique. However, the cost function is quadratic in the unknowns. The gradient of the cost function is therefore a linear function of the unknowns. The vanishing of its gradient must have a unique solution. Therefore, the solutions to both the minimization problems are identical.  
$\blacksquare$
\end{proof}

The theorem shows that the imposition of the zeros at the endpoints by the Lagrange multiplier gives nothing new when representing $u_{N}(x)$ as a finite sum in either the difference basis or the quadratic factor basis.


\subsection{Splitting the Least Squares Problem Into Two Via Parity}

An arbitrary function can always be split into its parts which are symmetric respect to reflection about the origin, $S(x)$, and antisymmetric with respect to reflection, $A(x)$ (Chapter 8 of \cite{Boyd99z}). Symmetry means $S(-x)=S(x), \forall x \in \Omega$, while $A(-x)= - A(x), \forall x\in \Omega$, where the $\Omega$ is the domain of a function. The parts are
$S(x) = (u(x) + u(-x))/2$ and $A(x) = (u(x) - u(-x))/2$.

If we apply this splitting to $u_{N}(x)$, the cost function becomes
\begin{eqnarray}
	\mathcal{J} & = & \mathcal{J}_{S} + \mathcal{J}_{A} + \lambda \Psi + \mu \Upsilon,
\end{eqnarray}
where
\begin{eqnarray*}
	\mathcal{J}_{S} &   =  &  < S - S_{N}, S - S_{N} >,  \quad 
	\mathcal{J}_{A}   =    < A - A_{N}, A - A_{N} >.
\end{eqnarray*}
where $\mathcal{J}_{S}$ is a function of the even degree spectral coefficients only while $\mathcal{J}_{A}$ is a function only of $\{ d_{1}, d_{3}, d_{5} \ldots \}$. After expanding the integrand of the original cost function to 
$<(S-S_{n}), (S-S_{n})> + <(A - A_{n}), (A - A_{n})> + <(S-S_{n}), (A - A_{n})> +<(A - A_{n}), (S-S_{n})>$, invoke the fact that the product of a symmetric function with an antisymmetric function is antisymmetric; the integral of an antisymmetric function over a symmetric interval is always zero.

The cost function is not completely decoupled because the constraints depend on both even and odd coefficients. However, both constraints are always zero at the solution. It follows that any linear combination of the constraints is also a legitimate constraint. Define
\begin{eqnarray}
	\Theta & = & (\Psi +\Upsilon)/2 = \sum_{n=0} d_{2n} ,\quad 
	\chi =  (\Psi - \Upsilon)/2 = \sum_{n=0} d_{2n+1}.
\end{eqnarray}

Least-squares is now split into two completely independent problems. One is to minimize, using only symmetric basis functions, 
\begin{eqnarray}
	\mathcal{J}_{S} + \lambda'    \,  \Theta
\end{eqnarray}
and the other, using only basis functions antisymmetric with respect to the origin, is to minimize
\begin{eqnarray}
	\mathcal{J}_{A} + \mu'   \,  \chi.
\end{eqnarray}
Since the methods of attack are similar for each, we shall only discuss the even parity problem in detail.

\subsection{Unconstrained Least Squares with  Quadratic-Factor Basis}

Define
\begin{eqnarray}	
	v^{quad,LS}_{N-2}(x) & = & \sum_{n=0}^{N-2}   c_{n}^{LS}  \, T_{n}(x)	,\\
	u^{quad,LS}_{N}(x) & =& \sum_{n=0}^{N-2}   c_{n}^{LS} \, (1-x^{2}) \, T_{n}(x)	,
\end{eqnarray}
and the cost function
\begin{eqnarray}	
	\mathfrak{J} = \frac{1}{2} < u(x) - u^{quad,LS}_{N}(x) , u(x) - u^{quad,LS}_{N}(x) >,
\end{eqnarray}
and the definition of function $v(x)$ is given in \eqref{Eqn:vdef}.
Then
\begin{eqnarray}	
	\mathfrak{J} & = & \frac{1}{2} < v(x) - v^{quad,LS}(x) , (1 - x^{2})^{2} \, v(x) - v^{quad,LS}(x) > \\	
	& = &   \frac{1}{2} \int_{-1}^{1}   \left(v(x) - v^{quad,LS}(x) \right)^{2} (1 - x^{2})^{3/2} \, dx	.
\end{eqnarray}
It follows that $v^{quad,LS}(x)$ is a standard polynomial approximation to $v(x)$, but the weight function is not the usual Chebyshev weight of $1/\sqrt{1 - x^{2}}$ but rather
$(1- x^{2})^{3/2}$. The orthogonal basis with this weight is the set of Gegenbauer polynomials of order 2. The Gegenbauer polynomials are defined as those polynomials with the orthogonality integral
\begin{equation}	
	\int_{-1}^{1}  (1-x^{2})^{m-1/2}  \hat{C}_{n}^{m}(x) \, \hat{C}_{k}^{m}(x)  dx = 0 ,\, \, \, k \neq n ,	
\end{equation}
where the subscript is the degree of the polynomial and here the polynomials are normalized so that $\hat{C}_{n}^{m}(1)=1$. (Warning: this is not the standard textbook normalization, but is convenient for comparing rates of convergence near the endpoints; we have added a caret to the symbol for the Gegenbauer polynomials to emphasize this.) 

The Gegenbauer coefficients are \emph{not} equal to the Chebyshev coefficients. However, Theorem 6 of \cite{BoydOP208}, which is a specialization of the theorem proved in $	a_{n}^{2} \sim \mathcal{O}( 1/ n^{2 \varphi-1}) $. Note that the $a_{n}^{2}$ are not the coefficients of $u(x)$ but rather are the coefficients of $v(x)=u(x)/(1-x^{2})$ which has branch points proportional to $(1-x^{2})^{\varphi-1}\, \ln( 1 - x^{2})$ instead of $(1-x^{2})^{\varphi}  \, \ln( 1 - x^{2})$ when $\varphi\in \mathbb{N}$; Theorem~\ref{Thm1} must be applied with $\varphi \rightarrow \varphi-1$ so that the coefficients of $v(x)$ decrease more slowly than those of $u(x)$ by a factor of $1/n^{2}$. The error near the endpoints is one order worse than the rate of convergence of the coefficients. For $\varphi\notin \mathbb{N}, \vartheta\in \mathbb{N}_{+}$,  by the Corollary 3.4 in \cite{Xiang21}, the coefficients decay rate of the standard Gegenbauer polynomials ($C_{n}^{m}(x)$ without the caret) expansion for the function $(1-x^{2})^{\varphi-1}\, \ln^{\vartheta}( 1 - x^{2})$ is  proportional to $\mathcal{O}(\ln^{\vartheta}(n)/n^{2\varphi +2m})$. By the equality $C_{n}^{m}(x) = \frac{\mathrm{\Gamma}(n+2m)}{\mathrm{\Gamma}(2m)\mathrm{\Gamma}(n+1)}\hat{C}_{n}^{m}(x)$, it is easy to see that the normalized Gegenbauer coefficients decay as
\begin{eqnarray}~\label{Eq3over2}  	
	a_{n}^{2} \sim \ln^{\vartheta}(n) / n^{2 \varphi-1},\quad n \gg 1.
\end{eqnarray}


The $N$-term truncation of the Gegenbauer series has an $L_{\infty}$ error norm for $v(x)$ of $\mathcal{O}(\ln^{\vartheta}(N)/N^{2 \varphi-2})$. Lemma~\ref{LemmaBD} allows us to convert the asymptotic behavior of the Gegenbauer coefficients into a bound on the slowness of the rate of convergence of the error norm.








%







\begin{theorem}~\label{ThLose1}
	Suppose that the coefficients $a_{n}^{m}$ of a spectral series in Gegenbauer polynomials $\hat{C}^{m}_{n}(x)$ or Chebyshev polynomials $T_{n}(x)$ ($m= 0 $) satisfy the bound
	\begin{eqnarray}	
		|a_{n}^{m} | \leq \,  W \, \dfrac{ \ln^{\vartheta} (n) }{n^{\kappa}}, \quad \forall n \geq 1, \quad \text{fixed} \,  \, m, \quad \vartheta \in \mathbb{N}\quad \kappa>1,
	\end{eqnarray}
	where $W$ is a positive constant, then the error in truncating the spectral series after the $N$-th term satisfies the inequality
	
	\begin{eqnarray}		
		\left|    v(x) - \sum_{n=0}^{N} a_{n}^{m} \, \hat{C}^{m}_{n}(x)   \right|   \leq \, W \, \dfrac{ \ln^{\vartheta}(N) }{ (\kappa-1)  \, N^{\kappa-1}   }	.	
	\end{eqnarray}
	
\end{theorem}

\begin{proof}
: By the Baszenski-Delvos Lemma \ref{LemmaBD}, the theorem is easy to be proved. $\blacksquare$
\end{proof}

The theorem (combined with the Tuan and Elliott's theorem \cite{TuanElliott72} for Gegenbauer coefficients) yields the maximum pointwise error for $v(x)$. The error for $u(x) \equiv (1- x^{2}) \, v(x)$ is

\begin{eqnarray}~\label{Old188}	
	E^{u}_{N}(x) =  \sum_{n=N+1}^{\infty} a^{2}_{n} \, (1 - x^{2}) \,\hat{C}^{2}_{n}(x)	.
\end{eqnarray}

To proceed further, we need two additional lemmas.

\begin{lemma}[Gegenbauer As Chebyshev Derivative]~\label{GegChebDeriv}
	Normalize the Gegenbauer polynomials so that each is one at the right endpoint. Then
	
	\begin{eqnarray}		
		\hat{C}_{n-k}^{k}(x) = \dfrac{1}{  			
			\prod_{j=0}^{k-1} \frac{n^2 - j^2}{2j+1}			
		} \dfrac{d^{k} T_{n}(x)}{dx^{k}}		.
	\end{eqnarray}
	
\end{lemma}

\begin{proof}
: It has long been known (18.9.19 on pg. 446 of \cite{NISTLibrary} )  that the $k$-th derivative of a Gegenbauer polynomial $\hat{C}_{n}^{m}(x)$ is proportional to $\hat{C}_{n-k}^{m+k}(x)$. It only remains to deduce the proportionality constant. Since the Gegenbauer polynomials are normalized to be one at the right endpoint, this constant must be the reciprocal of the value of the derivative at the origin which is known analytically to be (Appendix A of \cite{Boyd99z}).

\begin{equation}  	
	\left. \frac{d^k T_n(x)}{dx^k} \right|_{x= 1} = 	
	\prod_{j=0}^{k-1} \frac{n^2 - j^2}{2j+1} 	.
\end{equation}

$\blacksquare$
\end{proof}

Applying Lemma~\ref{GegChebDeriv}, the error (\ref{Old188}) in the variational approximation of $u(x)$ transforms to

\begin{eqnarray}~\label{hottie}	
	E^{u}_{N}(x) =  \sum_{n=N+1}^{\infty} a^{2}_{n} \, (1 - x^{2}) \,	
	\dfrac{3}{  (n+2)^{2} \left((n+2)^{2} - 1\right)} \,	
	\dfrac{d^{2} T_{n+2}(x)}{dx^{2}} 	.
\end{eqnarray}




Sergei N. Bernstein proved the following elegant theorem in a paper written in French \cite{Bernstein1913}. Here $P_{n}$ denotes the space of all polynomials whose degree is no more than $n$.

\begin{theorem}[Bernstein Polynomial Derivative Bound]
	If $P(x)$ is a polynomial of degree less than  or equal to $n$ in $P_{n}$, then for $k\le n$ and $x\in [-1,1]$,
	\begin{equation}\label{BersteinInequ}		
		\left| \dfrac{d^{k} P(x)}{dx^{k}} \right| \le B_{n}^{(k)}  \, || P(x) ||_{\infty},
	\end{equation}
where 
\begin{equation*}
	B_{n}^{(k)} = \sup_{P}\left\{ |P^{(k)}(x)| : ||P||_{\infty} \le 1\,\;  \text{and}\, \; P \in P_{n}\right\}.
\end{equation*}
Moreover,  when $n$ is large, for $x\in(-1,1)$,
\begin{eqnarray}\label{BernsPara}
	B_{n}^{(k)} \simeq \left( \dfrac{n}{\sqrt{1-x^{2}}} \right)^{k} , 		
	\qquad n \rightarrow \infty.
\end{eqnarray}	
\end{theorem}

In this theorem, the inequality holds with increasing precision in the asymptotic limit of increasing degree.	A complete proof in English is given by R. Whitley in \protect\cite{Whitley85}.

Multiplying the equation \eqref{BersteinInequ} by $(1-x^{2})$ and specializing to $k=2$ and $P(x)=T_{n}(x)$ gives

\begin{eqnarray}~\label{Txxbound}	
	\left| ( 1-x^{2} ) \, \dfrac{d^{2} T_{n}(x)}{dx^{2}} \right| \leq (1 + c) \, n^{2}	,
	\qquad n \rightarrow \infty,
\end{eqnarray}
where the small parameter satisfies $0<c<1$. It is easy to prove that 
\begin{equation}
	(1-x^2) \frac{d^{2} T_{n}(x)}{dx^{2}} = n \left(n\, T_{n}(x) - x\,U_{n-1}(x)\right).
\end{equation}
Thus, when $n\in \mathbb{N}$ and $n\ge 2$, it holds that
\begin{equation}
		\left| ( 1-x^{2} ) \, \dfrac{d^{2} T_{n}(x)}{dx^{2}} \right| \le 2\,n^{2}, \quad x\in[-1,1].
\end{equation}

\begin{theorem}[Error Bound for $u(x)$ in Least Squares/Quadratic-Factor Basis]~\label{Thebb}
	The error $E_{N}^{u}(x)$ in degree $N$ approximation in  the quadratic-factor basis using least squares with the inner product $< f(x), g(x) > = \int_{-1}^{1} f(x) \, g(x) dx$ satisfies the inequality
	
	\begin{eqnarray}		
		|E_{N}^{u}(x)| \leq \, 6 \, \dfrac{W}{2 \varphi}   \,\dfrac{\ln^{\vartheta}(N)}{ N^{2 \varphi}   }	.	
	\end{eqnarray}
	
\end{theorem}

\begin{proof}
: We have previously shown that, repeating (\ref{hottie}) here for clarity,

\begin{eqnarray}	
	E^{u}_{N}(x) =  \sum_{n=N+1}^{\infty} a^{2}_{n} \, \dfrac{3}{  (n+2)^{2} ((n+2)^{2} - 1)} \, (1 - x^{2}) \,	
	\dfrac{d^{2} T_{n+2}(x)}{dx^{2}}.  \qquad  \, \, \, \nonumber	
\end{eqnarray}
Recall that we previously demonstrated that $a_{n}^{2}$ are proportional to $\ln^{\vartheta}(n)/n^{2 \varphi -1}$ in (\ref{Eq3over2}). Applying the bound on the second derivative of the Chebyshev polynomials (\ref{Txxbound}), the error bound transforms to

\begin{eqnarray*}
	|E^{u}_{N}(x) |   &  \leq   & \, W  \sum_{n=N+1}^{\infty} \dfrac{\ln^{\vartheta}(n)} {n^{2 \varphi-1} } \,   \dfrac{6}{  (n+2)^{2}((n+2)^2-1)} \, (n+2)^{2} 
	 \leq   \, 6 \, W  \sum_{n=N+1}^{\infty} \dfrac{\ln^{\vartheta}(n)} {n^{2 \varphi-1+2} } .\,   	
\end{eqnarray*}
Applying Lemma~\ref{LemmaBD} with $k=2 \varphi +1 $ proves the theorem.  $\blacksquare$
\end{proof}

By \eqref{BernsPara}, it is not hard to see that when $N$ is large, one can obtain a sharper  estimate 
\begin{equation}
|E_{N}^{u}(x)| \le (3+c) \,\frac{W}{2 \varphi}   \,\frac{\ln^{\vartheta}(N)}{ N^{2 \varphi}}, \quad N\to \infty,
\end{equation}
where the small parameter is in the interval $(0,1)$.

\subsection{Least Squares with the Difference Basis}

In this basis, the square matrix $\mathbf{G}$ has elements
\begin{eqnarray} G_{mn} &   = & < (T_{2m}(x)-T_{2m-2}(x),  T_{2n}(x) - T_{2n-2}(x)  >, m=1, 2, \ldots; 
	n=1, 2,  \ldots \nonumber\\
	& = & \left\{  \begin{array}{ll}   \pi,\qquad\quad  & m=n >3 ,\\
     	3 \pi/2, &m=n=1, \\
		- \pi/2, &m=n+1, \\
		- \pi/2, & n=m+1. \\
	\end{array} \right. 
\end{eqnarray}
Thus, the $6 \times 6$ case is
\begin{eqnarray}  \frac{2}{\pi} \mathbf{G} =
	\begin{vmatrix}
		3 & -1 & 0 & 0 & 0 & 0 \\
		-1 & 2 & -1 & 0 & 0 & 0 \\
		0 & -1 & 2 & -1 & 0 & 0 \\
		0 & 0 & -1 & 2 & -1 & 0 \\
		0 & 0  & 0 & -1 & 2 & -1 \\
		0 & 0 & 0 & 0 & -1 & 2 
	\end{vmatrix}.
\end{eqnarray}
and with $a_{n}$ denoting Chebyshev coefficients of the usual infinite series, unconstrained to vanish at the endpoints,
\begin{eqnarray}  \frac{2}{\pi} \mathbf{f} =
	\begin{vmatrix}
		a_{2} - 2 a_{0}  \\
		a_{4} - a_{2}  \\
		a_{6} - a_{4}  \\
		a_{8} - a_{6}  \\
		a_{10} - a_{8} \\
		a_{12} - a_{10}  
	\end{vmatrix}.
\end{eqnarray}
Because of its sparsity, the matrix equation $\mathbf{G} \mathbf{d} = \mathbf{f}$, with the $b^{LS}_{2n}$ now denoting the elements of $\mathbf{d}$, can be written as the difference system
\begin{eqnarray}
	3 b^{LS}_{0} - b^{LS}_{2} & = &    a_{2} - 2 a_{0}, \\
	- b^{LS}_{2n-2} + 2 b^{LS}_{2n} - b^{LS}_{2n+2} & = &  a_{2n+2} - a_{2n}, \qquad n=1,2,\, \ldots \, , (N-2),\\
	- b^{LS}_{2N-2} + 2 b^{LS}_{2N} & = &  a_{2N} - a_{2N-2}.\, \, 
\end{eqnarray}
The solution is
\begin{eqnarray}
	b^{LS}_{2n} & = & \, - \,  \frac{1 - n/N}{1 + 1/(2N)}  \, \sum_{m=0}^{n} \,  a_{2m} 
	+ \frac{n + 1/2}{N+1/2} \sum_{m=n+1}^{N} a_{2m}, \qquad n=0,\ldots, (N-1) .
	\nonumber   
\end{eqnarray}

The infinite series limit, already analyzed in Sec.~\ref{SecRateCon}, is
\begin{eqnarray}
	\lim_{N \rightarrow \infty, \mbox{fixed} \, n} b_{2n} & = & \, - \,  \sum_{m=0}^{n} \,  a_{2m} 
	, \qquad n=0,\ldots, \infty .
	\nonumber   
\end{eqnarray}
If both $n$ and $N$ are large but finite, the solution simplifies to
\begin{eqnarray}
	b^{LS}_{2n} & = & \, - \,  \left(  1 - \frac{n}{N} \right)  \, \sum_{m=0}^{n} \,  a_{2m} 
	+ \frac{n}{N} \sum_{m=n+1}^{N} a_{2m}, \qquad n=0,\ldots, (N-1) .
	\nonumber   
\end{eqnarray}
Now  the Chebyshev coefficients of $u(x)$ must satisfy the condition $u(1)=0$ which demands
\begin{eqnarray}
	\sum_{m=0}^{n} \,  a_{2m}  = - \sum_{m=n+1}^{\infty} a_{2m}.
\end{eqnarray}
Similarly, the second sum in $b^{LS}_{2n}$ can be rewritten in terms of infinite summations as
\begin{eqnarray}
	\sum_{m=n+1}^{N} a_{2m} = \sum_{m=n+1}^{\infty} a_{2m} - \sum_{m=N+1}^{\infty} 
	a_{2m}.
\end{eqnarray}
Then 
\begin{eqnarray}
	b^{LS}_{2n} & = &  \,  \left(  1 - \frac{n}{N} \right)  \, \sum_{m=n+1}^{\infty} \,  a_{2m} 
	+ \frac{n}{N} \sum_{m=n+1}^{\infty} a_{2m} - \frac{n}{N} \sum_{m=N+1}^{\infty} \,
	a_{2m} ,
	\nonumber   \\
	b^{LS}_{2n} & = &  \,   \, \sum_{m=n+1}^{\infty} \,  a_{2m}  - \frac{n}{N} \sum_{m=N+1}^{\infty} \,
	a_{2m} .
	\nonumber   
\end{eqnarray}
Recall from Lemma~\ref{LemmaBD} that (\ref{EqBasz})
\begin{eqnarray}\label{EqBaszbis}
	\sum_{n=N+1}^{\infty} \frac{\ln^{\vartheta}(n)}{(2n)^{\kappa}}  \,  \sim \frac{1}{(\kappa - 1) N^{\kappa-1}}  \frac{ \ln^{\vartheta}(N) }{ 2^{\kappa} } \left\{ 1 + \mathcal{O} \left( \frac{1}{N}  \right) \right\}, \quad \vartheta \in\mathbb{N}.
\end{eqnarray}
If the $a_{n} \sim \mathcal{O}( \ln^{\vartheta}(n)/n^{\kappa} )$, then
\begin{eqnarray}
	b^{LS}_{2n} & \sim & \frac{A}{(\kappa-1)} \frac{ 1 }{2^\kappa} \left\{ \frac{ \ln^{\vartheta}(n) }{n^{\kappa-1}}  
	- \frac{n}{N} \frac{ \ln^{\vartheta}(N) }{ N^{\kappa-1} } 
	\right\} 
	\sim  \frac{A\ln^{\vartheta}(n) }{(\kappa-1) n^{\kappa-1}} \frac{1}{2^{\kappa}} \left\{ 1 - \frac{n^{\kappa}}{N^{\kappa}} \dfrac{\ln^{\vartheta}(N)}{\ln^{\vartheta}(n)} 
	\right\} .
\end{eqnarray}

The coefficients in the infinite series are $b^{LS}_{2n} \sim \mathcal{O}(\ln^{\vartheta}(n)/n^{\kappa-1})$, which is the same power law of 
rate of decay as for its least squares counterparts. However, the least-squares coefficients --- but not the infinite series coefficients ---
multiply the $\ln^{\vartheta}(n)/n^{\kappa}$ by $\left(1 - (\frac{n}{N})^{\kappa} \frac{ \ln^{\vartheta}(N) }{ \ln^{\vartheta}(n) }\right)$. On a log-log plot, the $b^{LS}_{2n}$ curve sharply downward as $n \rightarrow N$.

\subsection{Least Squares for Chebyshev Series with Lagrange Multiplier}

If a constraint is not built-in to the approximation $u_{N}(x)$, it can alternatively be added by means of a Lagrange multiplier. The goal is to enforce two boundary conditions, but a function can always be split by parity and then only one constraint for each symmetry is needed.

The goal of least squares is to minimize the ``cost function"
\begin{eqnarray}
	J = \frac{1}{2} < u(x) - S_{N}(x) , u(x) - S_{N}(x) > + \lambda \Psi,
\end{eqnarray}
where, for the even parity case,  
\begin{eqnarray}
	S_{N}(x)&=& \sum_{n=0}^{N} a^{LS}_{2n} T_{2n}(x),\quad 
	\Psi = \sum_{n=0}^{N} a^{LS}_{2n}.
\end{eqnarray}
Setting the gradients of the cost function with respect to all unknowns gives
\begin{eqnarray}
	\frac{\partial J}{\partial \lambda} = \Psi =0,
\end{eqnarray}
which merely insists that the constraint be satisfied, and also
\begin{eqnarray}
	\frac{\partial J}{\partial a^{LS}_{2 m}} = 0 = \lambda  - \, <u(x), T_{2m}(x)> + \sum_{n=0}^{N} a^{LS}_{2n} < T_{2 m}(x), T_{2n}(x) >.
\end{eqnarray}
Because of orthogonality of the Chebyshev polynomials and using the identities $<T_{0}(x), T_{0}(x)>=\pi$ and
$<T_{2n}(x), T_{2n}(x)>=\pi/2$ for $n \geq 1$, the equations simplify to 
\begin{eqnarray}
	\lambda & =&  \, <u(x), T_{0}(x)> -  a^{LS}_{0} \pi, \qquad 
	\lambda  =  \, <u(x), T_{2n}(x)> -  a^{LS}_{2n} \pi/2, \quad n \geq 1.
\end{eqnarray}
Let $a_{n}$ denote the Chebyshev coefficients of the infinite series for $u(x)$. Recall that $a_{0} = \frac{1}{\pi}<u(x), T_{0}(x)>$ and $a_{2n} = \frac{2}{\pi}<u(x), T_{2n}(x) > $. Then
\begin{eqnarray}
	\frac{1}{\pi} \lambda & =&  \, a_{0} -  a^{LS}_{0},  \qquad 
	\frac{2}{\pi}  \lambda  =  \, a_{2n} -  a^{LS}_{2n} , \quad n \geq 1.
\end{eqnarray}
Adding these equations and then invoking $\Psi=0$ gives
\begin{eqnarray}
	\lambda = \frac{\pi}{ (2 N + 1)} \, \sum_{n=0}^{N} a_{2n}.
\end{eqnarray}
The Chebyshev coefficients of the solution to the variational problem are then
\begin{eqnarray}
	a^{LS}_{0} & = & - \frac{1}{\pi} \,  \lambda + \, a_{0},   \qquad
	a^{LS}_{2n}  =  - \frac{2}{\pi}  \, \lambda +  \, a_{2n} , \qquad n \geq 1.
\end{eqnarray}

If the $a_{2n} \sim A \ln^{\vartheta} (2n) / (2n)^{\kappa} (\kappa>0, \vartheta\in\mathbb{N})$ as demanded by the Theorem \ref{Thm1}, then 
the error at the endpoints is
\begin{eqnarray}
	\Upsilon      \equiv    \sum_{n=0}^{N} a_{2n} 
	= - \sum_{n=N+1}^{\infty} a_{2n} 
	\approx    \mathcal{O}\left( \frac{\ln^{\vartheta}(N)}{N^{\kappa-1}} \right)
\end{eqnarray}
will be $\mathcal{O}\left(\ln^{\vartheta}(N)/N^{\kappa-1}\right)$, the same as the $L_{\infty}$ error norm of the 
Chebyshev series. (The error norm in fact \emph{is} $\Upsilon$ for some of our exemplary $u(x)$.)
It follows that
\begin{eqnarray}
	\lambda \sim \mathcal{O}(\ln^{\vartheta}(N)/N^{\kappa}), \quad \vartheta\in \mathbb{N}.
\end{eqnarray}

It is deserving to point out that the least squares approximation varies with the choice of weight function. On above the Chebyshev weight function is selected as $(1-x^2)^{-1/2}$ for all three bases. However, only the Chebyshev basis is orthogonal with the weight function, the two others are not. Next, the weight functions to make the difference basis and the quadratic basis orthogonal are respectively given in this section. 
\begin{theorem}
	If the weight function is chosen as $(1-x^2)^{-3/2}$, then the difference basis $\{\varsigma_n(x)\}$ are orthogonal, i.e.
	\begin{equation}
		\int_{-1}^{1}\varsigma_{m}(x)\varsigma_{n}(x)(1-x^2)^{-\frac{3}{2}}dx=2\pi \delta_{mn},\quad m,n=0,1,\cdots.
	\end{equation}
\end{theorem}
\begin{proof} :
\begin{equation*}
	\begin{aligned}
		\int_{-1}^{1} \varsigma_{n}(x)\varsigma_{m}(x)(1-x^2)^{-\frac{3}{2}} dx 
		&=\int_{0}^{\pi} \left[\cos((n+2)t)-\cos(nt)\right]\cdot\left[\cos((m+2)t)-\cos(mt)\right]\frac{1}{\sin^2(t)}dt\\
		&= 4\int_{0}^{\pi} \sin[(m+1)t]\sin(t)\cdot \sin[(n+1)t]\sin(t)\cdot \frac{1}{\sin^2(t)}dt\\
		&=4\int_{0}^{\pi} \sin[(m+1)t]\cdot \sin[(n+1)t]\cdot dt\\
		&=2\pi \delta_{mn}.
	\end{aligned}
\end{equation*}
$\blacksquare$
\end{proof}

Following the steps of least squares in Sec.\ref{Sec7_1}, like (\ref{EqGsolve}), for the difference basis with the weight function $(1-x^2)^{-3/2}$, one obtains
\begin{equation}
	\tilde{\mathbf{G}}\tilde{\mathbf{d}}=\tilde{\mathbf{f}},
\end{equation}
where 
\begin{equation*}
	\begin{aligned}
		\tilde{\bf{G}}_{mn}&=<\varsigma_{m}(x),\varsigma_{n}(x)>=\int_{-1}^{1} \varsigma_{n}(x)\varsigma_{m}(x)(1-x^2)^{-\frac{3}{2}} dx\\
		&=4\int_{0}^{\pi} \sin[(m+1)t]\cdot \sin[(n+1)t]\cdot dt\\
		&\approx 4\cdot \frac{\pi}{N_{col}} \sum_{k=0}^{N_{col}}\sin\left((m+1)\frac{(2k-1)\pi}{2N_{col}}\right)\sin\left((n+1)\frac{(2k-1)\pi}{N_{col}}\right)\\
		&=4\cdot \frac{\pi}{N_{col}}\cdot \frac{N_{col}}{2}\delta_{mn}=2\pi\delta_{mn},\quad m,n=0,1,\, \ldots\, ,N, \;N_{col}>N+1,\\
	\end{aligned}
\end{equation*}
which is a consequence of the orthogonality of the sine function with respect to the points $t_k=\frac{\pi (2k-1)}{N_{col}}, k=1,2,\cdots, N_{col}$,
\begin{equation*}
	\sum_{k=1}^{N_{col}}\sin\left((m+1)\frac{(2k-1)\pi}{2N_{col}}\right)\sin\left((n+1)\frac{(2k-1)\pi}{2N_{col}}\right)=\frac{N_{col}}{2}\delta_{mn},
\end{equation*}
and 
\begin{equation*}
	\begin{aligned}
		\mathbf{\tilde{f}_n}=<u(x),\varsigma_n(x) >&=\int_{-1}^{1}\varsigma_{n}(x)u(x)(1-x^2)^{-\frac{3}{2}}dx\\
		&=\int_{0}^{\pi} \varsigma_{n}(\cos(t))u(\cos(t))\cdot\frac{1}{\sin^2(t)}dt\\	
		&\approx \frac{\pi}{N_{col}}\sum_{k=1}^{Ncol}\varsigma_{n}\left(\cos(\frac{2k-1}{2N_{col}}\pi)\right)u\left(\cos(\frac{2k-1}{2N_{col}}\pi)\right)\cdot\sin^{-2}\left(\frac{2k-1}{2N_{col}}\pi\right).
	\end{aligned}
\end{equation*}

When the number of interpolation $N_{col}$ is more than the number of basis $N+1$, the coefficients of difference basis decrease as $\mathcal{O}(\ln^{\vartheta}(n)/n^{2\varphi})$ as $n\to \infty$, which obeys the same law of the counterpart coefficients in infinite series truncation as is shown in Fig.\ref{FigOP279_Coeffs_InfiniteSeries}. There is no curl up or curl down as $n\rightarrow N$. Thus the error norm is also the same as the error norm of the infinite series truncation.  

In fact, to approximate the function $u(x)$, using the difference basis $\varsigma_{n}(x)$ with weight function $(1-x^2)^{-3/2}$ is equivalent to using the second Chebyshev function $U_n(x)$ with the weight function $\sqrt{1-x^2}$.

\begin{theorem}
	If the weight function is chosen as $(1-x^2)^{-5/2}$, then the quadratic basis $\{\varrho_n(x)\}$ are orthogonal, i.e.
	\begin{eqnarray}
		\int_{-1}^{1}\varrho_{m}(x)\varrho_{n}(x)(1-x^2)^{-\frac{5}{2}}dx=
		\left\{
		\begin{aligned}
			& \pi ,\quad &m=n=0,\\
			&\frac{1}{2}\pi \delta_{mn},\quad &m,n\in\mathbb{N}_{+}.			
		\end{aligned}	
		\right.
	\end{eqnarray}
\end{theorem}

The theorem is easy to be proved. Similar to the procedure of least squares for difference basis with weight function $(1-x^2)^{-3/2}$, one can also conclude that the least squares coefficients for quadratic basis with weight function $(1-x^2)^{-5/2}$ decrease as $\mathcal{O}( A(n)/n^{2\varphi-1} )$ as $n\to\infty$. The error decreases as $\mathcal{O} ( |A(N)| /N^{2\varphi-1} )$, which is also same as the error of the infinite series truncation for the same basis.  In the rest of this section, we still use the inner product $<\cdot,\cdot>$ mentioned in Sec. \ref{Sec7_1}.

\section{Comparing Different Approximations Using the Difference Basis}

The spectral coefficients and error norms are so similar that the most illuminating way to compare 
them is to tabulate ratios. Table~\ref{Tabbnratio} shows that when $n \ll N$, 
$b_{n} \approx b_{n}^{I} \approx    b_{n}^{LS}$. When $n$ nears $N$, the interpolation 
coefficients swell to nearly double those of the infinite series while $b_{LS} \ll b_{n}$.

\begin{table}[H]
	\vspace{-0.1in}
	\caption{~\label{Tabbnratio}Coefficient ratios for the difference basis, $\varsigma_{n}=T_{n+2}(x) - T_{n}(x)$, for $(1+x/2) (1-x^2)^{\varphi}\ln(1-x^2)^{\vartheta}$ for
		$\vartheta=1$ and $\varphi=2$. The $b_{n}$ are the coefficients in the infinite series, $b_{n}^{I}$ are the coefficients of the interpolant using 100 collocation points and $b_{n}^{LS}$ are the result of least squares with integration as the inner product.}
		\begin{center} 
		\begin{tabular}{|l|llllllllllllll|}
			\hline
			$n$    &10 & 20 & 30 &40 &50 &60 &70 &80  &90  & 92  &94  & 96  &98 &99 \\ 
			\hline 
			$b_{n}^{I}/b_{n}$ &1.00 & 1.00 &1.00&1.00&1.01&1.03&1.08&1.19&1.43&1.51&1.60&1.71&1.83&1.90\\
			\hline
			$b_{n}^{LS}/b_{n}$ &1.00&1.00&1.00&0.99&0.97&0.93&0.85&0.71&0.47&0.40&0.34&0.26&0.18&0.095\\
			\hline
		\end{tabular} 
	\end{center}
\end{table}
\vspace{-0.1in}
\begin{table}[h]
	\caption{~\label{Tabdifferrratio} Same as previous table except that the ratios 
		are now  of errors in the $L_{\infty}$ norm, and these are  listed versus the truncation $N$ rather than degree $n$.}
	\begin{center} 
		\begin{tabular}{|l|llllllllll|}
			\hline
			$N$ &10  & 20 & 30 & 40 & 50 & 60 & 70 & 80 & 90 &100 \\
			\hline
			$E_{N}^{interp}/E_{N}$ &1.98 & 1.96 & 1.96 & 1.97 & 1.94 & 1.93 & 1.98 & 1.98 & 1.96 & 1.85\\
			\hline
			$E^{LS}_{N}/E_{N}$ & 0.96 & 1.03 & 1.06 & 1.07 & 1.08 & 1.09 & 1.07 & 1.09 & 1.04 & 1.07\\
			\hline
		\end{tabular} 
	\end{center}
\end{table}

Table~\ref{Tabdifferrratio} compares the ratio of error norms. Least squares with integration as the inner product is only slightly worse than the truncation of the infinite series ( less than 10 \%). The maximum pointwise error for interpolation is roughly double that of truncation of the infinite series, independent of $N$.

\section{Comparing Different Quadratic-Factor Basis Approximations}

\begin{figure}
	\centerline{\includegraphics[scale=0.65]{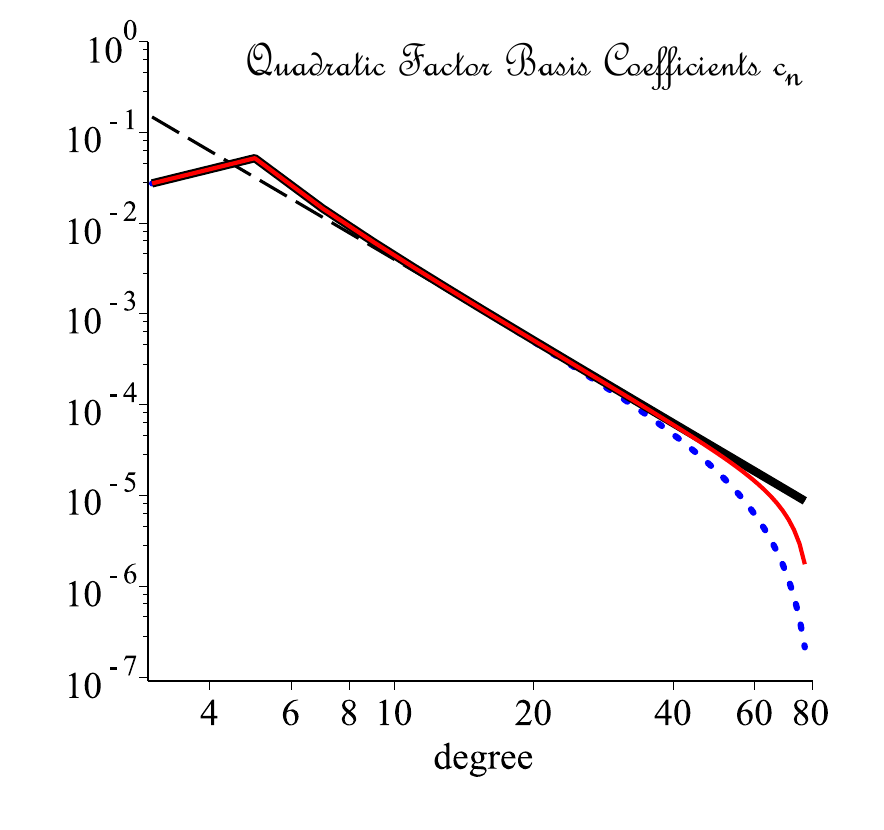}}
	\caption{Odd degree coefficients versus degree for approximations using the quadratic-factor basis for 
		$u(x)=(1+x/2)  (1- x^{2})^{2} \, \ln(1-x^{2})$, the same function as employed in the previous figure, $N=80$.  The thick black curve is the coefficients $c_{n}$  of the infinite series. The thin red curve connects the absolute values of the coefficients of the 79 point interpolant in the quadratic-factor basis. The blue dotted curve is the coefficients of least-squares with  integral inner product.  The black dashed 
		line is  $4/n^{3}$, proportional to $n^{\kappa-2}$.
	}
	\label{FigCoompare_three_cn_quadfactor_basis_COEFFS}
\end{figure}

Fig.~\ref{FigCoompare_three_cn_quadfactor_basis_COEFFS} shows that the coefficients of 
all three approximation schemes in the basis $\varrho_{n}(x) = (1 - x^{2}) \, T_{n}(x)$
has the same slope, $1/n^{2 \varphi - 1}$ over most of the range in degree. The interpolant's coefficients and those obtained by least-squares with the inner product $< f(x), g(x) >= \int_{-1}^{1}  f(x) \, g(x)/\sqrt{1-x^{2}} \,dx$ both bend sharply downward as $n \rightarrow N$.

\begin{figure}
	\centerline{\includegraphics[scale=0.65]{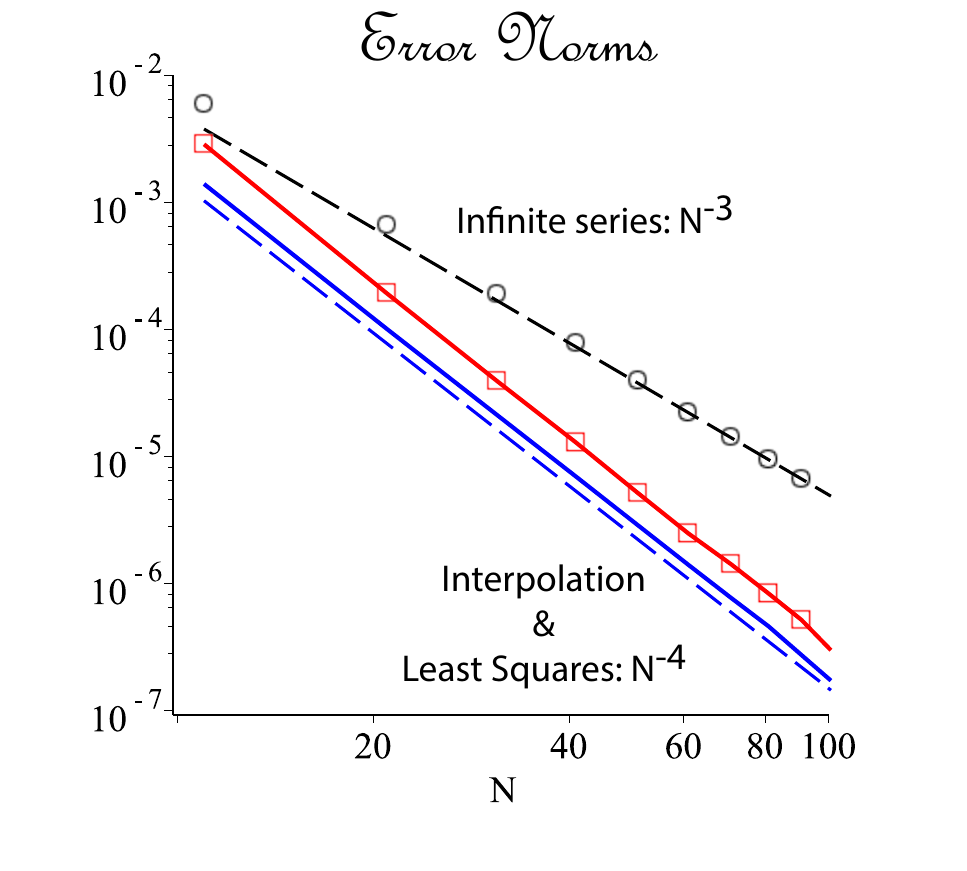}}
	\caption{Quadratic-factor basis, same as the previous figure but showing error norms versus $N$ instead of coefficients versus $n$. Black solid curve: errors in truncation of the infinite series in the basis $\varrho_{n}(x)$.}
	\label{FigCoompare_three_cn_quadfactor_basis_CERRORNORMS}
\end{figure}

How do these fast-tail decreases affect the error norms?
Fig.~\ref{FigCoompare_three_cn_quadfactor_basis_CERRORNORMS} provides an answer.

Aliasing error, which produces the downward curve in the spectral coefficients for 
interpolation in Fig.~\ref{FigCoompare_three_cn_quadfactor_basis_COEFFS}, is generally 
regarded as a Bad Thing.  Therefore, the even sharper deviation from a power law for the least squares coefficients should be an even Badder Thing. Actually, the error norms associated 
with the downward curving spectral coefficients decrease faster by $O(N)$ than the error norm of the truncated infinite 
series with its pure power law (black straight line in Fig.~\ref{FigCoompare_three_cn_quadfactor_basis_COEFFS}).

We have no explanation. However, note that some acceleration methods such as Euler acceleration  \cite{BoydMoore86,MorseFeshbach53,Pearce78}
taper the high degree coefficients to improve accuracy.  Something similar seems to be happening with aliased spectral series.

%

\section{Summary}

The concern of this paper has been to address the Chebyshev expansion of the weak singularity functions on three bases, both theoretically and computationally. The main results are concluded in the following.

\begin{enumerate}
	\item The coefficients and errors of several kinds of approximations are summarized in Table \ref{Tab3}.
	
	\begin{table}[H]
		\caption{Results for $(1-x^2)^{\varphi}\ln^{\vartheta}(1-x^2)$ with $\varphi>0$ and $\vartheta\in\mathbb{N}$. The labels ``u" or "d" denote that the coefficients $a_{n}$  curl up or 
			curl down as $n   \rightarrow N$, deviating from the correct asymptotic line because of 
			aliasing errors as described in Theorem~\ref{Th6}. The expression of $A(n)$ is given in Theorem 
		\ref{Thm1}. TS, IT, LS, and B.C.s represent Truncated Series, Interpolation, Least Squares and Boundary Conditions.}\label{Tab3}
		\begin{center}
		{\footnotesize
			\hspace*{-0.05in}~\begin{tabular}{|c|c|c|c|c|}
				\hline
				\linespread{1.2}
				Bases & Chebyshev & Difference & Quadratic & Chebyshev Lagrange \\  \hline
				TS: Coeffs & $a_{n} \sim A(n)/n^{2\varphi+1}$
				&   $b_{n} \sim A(n)/ n^{2\varphi}$ & $c_{n} \sim- A(n)/n^{2\varphi-1}$ &- \\
				
				TS : Error  & $|A(N)|/N^{2\varphi}$ & $ |A(N)|/N^{2\varphi}$ & $|A(N)|/N^{2\varphi-1}$ & - \\  \hline

				IT: Coeffs & $a^{I}_{n} \sim A(n)/n^{2\varphi+1}$ (d) & $b^{I}_{n} \sim A(n)/ n^{2\varphi} $ (u)   & $c^{I}_{n} \sim -A(n)/n^{2\varphi-1}$  (d)& - \\
				
				IT: Error & $|A(N)|/N^{2\varphi}$ &  $|A(N)|/N^{2\varphi}$   & $|A(N)|/N^{2\varphi}$ & - \\    \hline

				LS : Coeffs & $a_{n}^{LS} \sim  A(n)/n^{2\varphi+1}$ & $b_{n}^{LS} \sim  A(n)/n^{2\varphi}$ (d) & $c_{n}^{LS} \sim - A(n)/n^{2\varphi-1}$ (d) & $a_{n}^{CL} \sim  A(n)/n^{2\varphi+1}$ \\
				
				LS: Error  & $|A(N)|/N^{2\varphi}$ &  $|A(N)|/N^{2\varphi}$   & $|A(N)|/N^{2\varphi}$ & $|A(N)|/N^{2\varphi}$ \\    \hline

				B.C.s & Not imposed & Satisfied & Satisfied & Imposed by \\
				&    &   &   &  Lagrange Multiplier   \\                
				
				\hline
		\end{tabular} }
	\end{center}
\end{table} 

	\item There are \emph{two} distinct interpolants on the \emph{roots} grid, but the interpolant 
	on the \emph{Lobatto} (endpoint-including) grid is always \emph{unique}.

	\item The error norms in $N$-point interpolation on the roots grid are \emph{identical} for 
	all three basis sets, that is,
	\begin{eqnarray}
		{\rm Chebyshev}{\it } =
		\mbox{difference basis} =
		\mbox{quadratic-factor basis} \sim \mathcal{O}\left(|A(N)|/N^{2 \varphi}\right).
	\end{eqnarray}
	
	\item The pointwise errors for interpolation using the quadratic-factor basis and difference basis are identical for all $x$ because $u_{N}^{diff}(x)=u_{N}^{quad}(x)$ for all $x$.
	
	\item The pointwise error in standard Chebyshev interpolation, unconstrained by $u_{N}(\pm 1)=0$, is \emph{different} from the errors (not error norms) of the constrained basis sets, the quadratic-factored basis and the difference basis; the errors of the constrained basis sets are \emph{nearly-uniform} in $x$ whereas the Chebyshev error is one order smaller than that of the constrained bases except in narrow \emph{boundary layers} where the Chebyshev error rises to equal that of the constrained bases.
	


	\item If the Chebyshev coefficients decay as $a_{n} \sim A\ln^{\vartheta}(n)/n^{\kappa}$ where $A$ is a constant, $\vartheta$ is a nonnegative integer and $\kappa>0$, then
	
	(a). For small degree $1<n\ll N$,  the relative error in the Chebyshev coefficient is
	\begin{eqnarray}
		\frac{ | \mathcal{E}_{n}| }{|a_{n}|}  &  \leq  &   \frac{1}{2^{\kappa-1} }  \frac{  n^{\kappa} }{  N^{\kappa}  } \frac{\ln^{\vartheta}(N)}{\ln^{\vartheta}(n)} .
	\end{eqnarray}
	
	(b). For $n=N-m$ when $m$ is small,
	 the relative error is
	\begin{eqnarray}
		\frac{ | \mathcal{E}_{N-m}| }{|a_{N-m}|} \sim 1 + \mathcal{O}\left( \frac{\kappa m}{N} \right).
	\end{eqnarray}
	When $\vartheta=0$, the $a^{I}_{n}$ is plotted versus $n$ with logarithmic axes, the curve on the log-log plot should, for power-law decay, asymptote to a straight line. The aliasing errors create a sharp downward turn in $a_{n}^{I} $ as $n \rightarrow N$. (The coefficients for the difference basis exhibit a sharp upturn for similar reasons.)
	
	\item The value of the derivatives at the endpoints is $\mathcal{O}(N^{2})$ for Chebyshev polynomials and for the quadratic-factor basis, but only $\mathcal{O}(N)$ for the difference basis.
\end{enumerate}

The most important conclusion is that those different choices of approximation schemes and bases can alter the rate of convergence by a factor of $N$ or $N^{2}$. For series that converge
proportionally to small inverse powers of $N$ due to weak endpoint singularities, this is significant. Knowing how to solve the singular problems using spectral methods is important, but giving which basis is optimal seems more practically significant, especially in high-dimensional spaces.

\section*{Acknowledgments}
This work was supported by the National Science Foundation of the U. S. under DMS-1521158 , the National Natural Science Foundation of China (No. 12101229), the Hunan Provincial Natural Science Foundation of China (No.2021JJ40331), and the Chinese Scholarship Council 201606060017, 202106720024. We would thank the four anonymous referees, whose comments greatly improved the paper.

	\bibliographystyle{siam}  

\bibliography{OP279_Corner_Singularity}
%
%

\end{document}